\documentclass[11pt]{amsart}
\usepackage{amsmath,amssymb,mathrsfs}
\newtheorem{theorem}{Theorem}[section]
\newtheorem{proposition}[theorem]{Proposition}
\newtheorem{corollary}[theorem]{Corollary}
\newtheorem{lemma}[theorem]{Lemma}
\newtheorem{conjecture}[theorem]{Conjecture}
\theoremstyle{definition}
\newtheorem{definition}[theorem]{Definition}
\begin{document}

\title[Minimal surfaces in $S^3$]{Minimal surfaces in $S^3$: a survey of recent results}
\author{Simon Brendle}
\dedicatory{Dedicated to Professor Blaine Lawson on the occasion of his seventieth birthday.}
\address{Department of Mathematics \\ Stanford University \\ Stanford, CA 94305}
\begin{abstract}
In this survey, we discuss various aspects of the minimal surface equation in the three-sphere $S^3$. After recalling the basic definitions, we describe a family of immersed minimal tori with rotational symmetry. We then review the known examples of embedded minimal surfaces in $S^3$. Besides the equator and the Clifford torus, these include the Lawson and Kapouleas-Yang examples, as well as a new family of examples found recently by Choe and Soret. We next discuss uniqueness theorems for minimal surfaces in $S^3$, such as the work of Almgren on the genus $0$ case, and our recent solution of Lawson's conjecture for embedded minimal surfaces of genus $1$. More generally, we show that any minimal surface of genus $1$ which is Alexandrov immersed must be rotationally symmetric. We also discuss Urbano's estimate for the Morse index of an embedded minimal surface and give an outline of the recent proof of the Willmore conjecture by Marques and Neves. Finally, we describe estimates for the first eigenvalue of the Laplacian on a minimal surface. 
\end{abstract}
\thanks{The author was supported in part by the National Science Foundation under grant DMS-1201924. Part of this work was carried out while the author visited Cambridge University. He is very grateful to the Department of Mathematics and Mathematical Statistics at Cambridge University for its hospitality, and to Professors Robert Kusner and Neshan Wickramasekera for helpful discussions.}

\maketitle 

\section{Introduction}

Minimal surfaces are among the most important objects studied in differential geometry. Of particular interest are minimal surfaces in manifolds of constant curvature, such as the Euclidean space $\mathbb{R}^3$, the hyperbolic space $\mathbb{H}^3$, and the sphere $S^3$. The case of minimal surfaces in $\mathbb{R}^3$ is a classical subject; see e.g. \cite{Lawson-book} for an introduction. In this paper, we will focus on the case when the ambient space is the sphere. Throughout this paper, we will identify $S^3$ with the unit sphere in $\mathbb{R}^4$; that is, 
\[S^3 = \{x \in \mathbb{R}^4: x_1^2+x_2^2+x_3^2+x_4^2=1\}.\] 
Let $\Sigma$ be a two-dimensional surface in $S^3$, and let $\nu$ be a unit normal vector field along $\Sigma$. In other words, we require that $\nu$ is tangential to $S^3$, but orthogonal to the tangent plane to $\Sigma$. The extrinsic curvature of $\Sigma$ is described by a symmetric two-tensor $h$, which is referred to as the second fundamental form of $\Sigma$. The second fundamental form is defined by 
\[h(e_i,e_j) = \langle D_{e_i} \nu,e_j \rangle,\] 
where $\{e_1,e_2\}$ is an orthonormal basis of tangent vectors to $\Sigma$. The eigenvalues of $h$ are referred to as the principal curvatures of $\Sigma$. The product of the principal curvatures depends only on the intrinsic geometry of $\Sigma$; in fact, the Gauss equations imply that 
\[K=1+\lambda_1\lambda_2,\] 
where $\lambda_1,\lambda_2$ are the principal curvatures of $\Sigma$ and $K$ denotes the intrinsic Gaussian curvature. Moreover, the sum of the principal curvatures is referred to as the mean curvature of $\Sigma$:  
\[H = \lambda_1+\lambda_2 = \sum_{i=1}^2 \langle D_{e_i} \nu,e_i \rangle.\] 
Geometrically, the mean curvature can be viewed as an $L^2$-gradient of the area functional; more precisely, given any smooth function $u$ on $\Sigma$, we have 
\[\frac{d}{dt} \text{\rm area}(\Sigma_t) \Big |_{t=0} = \int_\Sigma H \, u,\] 
where 
\[\Sigma_t = \{\cos(t \, u(x)) \, x + \sin(t \, u(x)) \, \nu(x): x \in \Sigma\}.\] 
This motivates the following definition: 

\begin{definition} 
A two-dimensional surface $\Sigma$ in $S^3$ is said to be a minimal surface if the mean curvature of $\Sigma$ vanishes identically.
\end{definition} 

The condition that $\Sigma$ is minimal can be rephrased in several equivalent ways:

\begin{theorem} 
Let $\Sigma$ be a two-dimensional surface in $S^3$. Then the following statements are equivalent: 
\begin{itemize}
\item $\Sigma$ is a minimal surface. 
\item $\Sigma$ is a critical point of the area functional. 
\item The restrictions of the coordinate functions in $\mathbb{R}^4$ are eigenfunctions of the operator $-\Delta_\Sigma$ with eigenvalue $2$; that is, $\Delta_\Sigma x_i + 2x_i = 0$ for $i \in \{1,2,3,4\}$.
\end{itemize}
\end{theorem}

In the following, we will be interested in closed minimal surfaces. While there are no closed minimal surfaces in $\mathbb{R}^3$, there do exist interesting examples of closed minimal surfaces in $S^3$. The simplest example of such a surface is the equator, which is defined by
\[\Sigma = \{x \in S^3 \subset \mathbb{R}^4: x_4=0\}.\] 
The principal curvatures of this surface are both equal to zero. In particular, the resulting surface $\Sigma$ is minimal. Moreover, the equator has constant Gaussian curvature $1$, and $\Sigma$ equipped with its induced metric is isometric to the standard sphere $S^2$.

Another basic example of a minimal surface in $S^3$ is the so-called Clifford torus. This surface is defined by 
\[\Sigma = \Big \{ x \in S^3 \subset \mathbb{R}^4: x_1^2+x_2^2=x_3^2+x_4^2=\frac{1}{2} \Big \}.\] 
In this case, the principal curvatures are $1$ and $-1$, so the mean curvature is again equal to zero. Moreover, the intrinsic Gaussian curvature vanishes, and $\Sigma$ equipped with its induced metric is isometric to the flat torus $S^1(\frac{1}{\sqrt{2}}) \times S^1(\frac{1}{\sqrt{2}})$.

In the 1960s, Lawson \cite{Lawson1} constructed an infinite family of immersed minimal tori in $S^3$ which fail to be embedded (see also \cite{Hsiang-Lawson}). Moreover, immersed minimal tori in $S^3$ have been studied intensively using integrable systems techniques; see e.g. \cite{Bobenko} or \cite{Hitchin}. In the remainder of this section, we describe a family of immersed minimal tori in $S^3$ which are rotationally symmetric. These surfaces are not embedded, but they turn out to be immersed in the sense of Alexandrov.

\begin{definition} 
\label{definition.of.alexandrov.immersion}
A map $F: \Sigma \to S^3$ is said to be Alexandrov immersed if there exists a compact manifold $N$ and an immersion $\bar{F}: N \to S^3$ such that $\Sigma = \partial N$ and $\bar{F}|_\Sigma = F$.
\end{definition}

The notion of an Alexandrov immersion was introduced by Alexandrov \cite{Alexandrov} in connection with the study of constant mean curvature surfaces in Euclidean space, and has since been studied by many authors; see e.g. \cite{Korevaar-Kusner-Ratzkin}, \cite{Korevaar-Kusner-Solomon}, \cite{Kusner-Mazzeo-Pollack}. Using the method of moving planes, Alexandrov \cite{Alexandrov} was able to show that any closed hypersurface in Euclidean space which has constant mean curvature and is Alexandrov immersed must be a round sphere. 

The following result was pointed out to us by Robert Kusner:

\begin{theorem} 
\label{alexandrov.immersed.tori}
There exists an infinite family of minimal tori in $S^3$ which are Alexandrov immersed, but fail to be embedded. 
\end{theorem}

\textbf{Proof of Theorem \ref{alexandrov.immersed.tori}.} 
We consider an immersion of the form 
\[F(s,t) = (\sqrt{1-r(t)^2} \, \cos s,\sqrt{1-r(t)^2} \, \sin s,r(t) \, \cos t,r(t) \, \sin t),\] 
where $r(t)$ is a smooth function which takes values in the interval $(0,1)$. Clearly, 
\[g \Big ( \frac{\partial}{\partial s},\frac{\partial}{\partial s} \Big ) = 1-r(t)^2\] 
and 
\[g \Big ( \frac{\partial}{\partial t},\frac{\partial}{\partial t} \Big ) = \frac{r'(t)^2+r(t)^2 \, (1-r(t)^2)}{1-r(t)^2}.\]  
Moreover, the unit normal vector field to the surface is given by 
\begin{align*} 
\nu(s,t) 
&= \frac{r(t)^2 \, \sqrt{1-r(t)^2}}{\sqrt{r'(t)^2+r(t)^2 \, (1-r(t)^2)}} \, (\cos s,\sin s,0,0) \\ 
&- \frac{r(t) \, (1-r(t)^2)}{\sqrt{r'(t)^2+r(t)^2 \, (1-r(t)^2)}} \, (0,0,\cos t,\sin t) \\ 
&- \frac{r'(t)}{\sqrt{r'(t)^2+r(t)^2 \, (1-r(t)^2)}} \, (0,0,\sin t,-\cos t). 
\end{align*} 
Hence, the second fundamental form satisfies 
\[h \Big ( \frac{\partial}{\partial s},\frac{\partial}{\partial s} \Big ) = \frac{r(t)^2 \, (1-r(t)^2)}{\sqrt{r'(t)^2+r(t)^2 \, (1-r(t)^2)}}\] 
and 
\[h \Big ( \frac{\partial}{\partial t},\frac{\partial}{\partial t} \Big ) = \frac{r(t) \, (1-r(t)^2) \, r''(t) - (2-3 \, r(t)^2) \, r'(t)^2 - r(t)^2 \, (1-r(t)^2)^2}{(1-r(t)^2) \, \sqrt{r'(t)^2+r(t)^2 \, (1-r(t)^2)}}.\] 
Therefore, the mean curvature vanishes if and only if $r(t)$ satisfies the differential equation 
\begin{equation} 
\label{ode}
(1-r(t)^2) \, r(t) \, r''(t) = (1-2 \, r(t)^2) \, (2 \, r'(t)^2+r(t)^2 \, (1-r(t)^2)). 
\end{equation}
The equation (\ref{ode}) implies that 
\[\frac{d}{dt} \Big ( \frac{r'(t)^2}{r(t)^4 \, (1-r(t)^2)^2} + \frac{1}{r(t)^2 \, (1-r(t)^2)} \Big ) = 0,\] 
hence 
\begin{equation} 
\label{conservation.law}
\frac{r'(t)^2}{r(t)^4 \, (1-r(t)^2)^2} + \frac{1}{r(t)^2 \, (1-r(t)^2)} = \frac{4}{c^2}
\end{equation}
for some constant $c$. This conserved quantity can also be obtained in a geometric way via Noether's principle; to that end, one applies the formula for the first variation of area to the ambient rotation vector field $K = (0,0,x_4,-x_3)$. 

The function $r(t) = \frac{1}{\sqrt{2}}$ is an equilibrium solution of (\ref{ode}), and the corresponding minimal surface is the Clifford torus. In view of (\ref{conservation.law}), any nearby solution of the differential equation (\ref{ode}) is periodic. Moreover, the period is given by the formula
\[T(c) = 2 \int_{\underline{x}(c)}^{\overline{x}(c)} \frac{c}{x \, \sqrt{1-x^2} \, \sqrt{4x^2(1-x^2) - c^2}} \, dx\] 
where $c<1$. Here, $\underline{x}(c)$ and $\overline{x}(c)$ are defined by 
\[\underline{x}(c) = \sqrt{\frac{1-\sqrt{1-c^2}}{2}}\] 
and 
\[\overline{x}(c) = \sqrt{\frac{1+\sqrt{1-c^2}}{2}}.\] 
Note that $T(c) \to \sqrt{2} \, \pi$ as $c \nearrow 1$.

We now choose the parameter $c<1$ in such a way that the ratio $\frac{2\pi}{T(c)}$ is rational. This implies that we can find a positive integer $k$ such that $\frac{2\pi k}{T(c)}$ is an integer. As a result, we obtain a solution $r(t)$ of the differential equation (\ref{ode}) satisfying $r(t + 2\pi k) = r(t)$. Having chosen $c$ and $r(t)$ in this way, the map 
\begin{align*} 
&F: [0,2\pi] \times [0,2\pi k] \to S^3, \\ 
&\qquad (s,t) \mapsto (\sqrt{1-r(t)^2} \, \cos s,\sqrt{1-r(t)^2} \, \sin s,r(t) \, \cos t,r(t) \, \sin t) 
\end{align*}
defines a minimal immersion of the torus $S^1 \times S^1$ into $S^3$. 

It remains to show that $F$ is an Alexandrov immersion. To see this, we consider the map 
\begin{align*} 
&\bar{F}: B^2 \times [0,2\pi k] \to S^3, \\ 
&\qquad (\xi,t) \mapsto \frac{(\sqrt{1-r(t)^2} \, \xi_1,\sqrt{1-r(t)^2} \, \xi_2,r(t) \, \cos t,r(t) \, \sin t)}{\sqrt{(1-r(t)^2) \, |\xi|^2 + r(t)^2}}, 
\end{align*}
where $B^2 = \{\xi \in \mathbb{R}^2: |\xi| \leq 1\}$. Since $r(t)$ is periodic with period $2\pi k$, the map $\bar{F}$ defines an immersion of the solid torus $B^2 \times S^1$ into $S^3$. Since $\bar{F}(\cos s,\sin s,t) = F(s,t)$, the map $F$ is an Alexandrov immersion. \\

\section{Examples of embedded minimal surfaces in $S^3$}

While Theorem \ref{alexandrov.immersed.tori} provides a large family of Alexandrov immersed minimal surfaces in $S^3$, it is a difficult problem to construct examples of minimal surfaces which are embedded. In fact, for a long time the equator and the Clifford torus were the only known examples of embedded minimal surfaces in $S^3$. This changed dramatically in the late 1960s, when Lawson discovered an infinite family of embedded minimal surfaces of higher genus: 

\begin{theorem}[H.B.~Lawson, Jr. \cite{Lawson2}]
\label{lawson.surfaces}
Given any pair of positive integers $m$ and $k$, there exists an embedded minimal surface $\Sigma$ in $S^3$ of genus $mk$. In particular, there exists at least one embedded minimal surface of any given genus $g$, and there are at least two such surfaces unless $g$ is a prime number.
\end{theorem} 

\textbf{Sketch of the proof of Theorem \ref{lawson.surfaces}.} 
For $i \in \mathbb{Z}_{2(k+1)}$ and $j \in \mathbb{Z}_{2(m+1)}$, we define  
\[P_i = \Big ( \cos \frac{\pi i}{k+1},\sin \frac{\pi i}{k+1},0,0 \Big )\] 
and 
\[Q_j = \Big ( 0,0,\cos \frac{\pi j}{m+1},\sin \frac{\pi j}{m+1} \Big ).\] 
Moreover, let $A = \mathbb{Z}_{2(k+1)} \times \mathbb{Z}_{2(m+1)}$ and 
\begin{align*} 
A_{\text{\rm even}} 
&= \{(i,j) \in \mathbb{Z}_{2(k+1)} \times \mathbb{Z}_{2(m+1)}: \text{\rm $i$ and $j$ are both even}\} \\ 
&\cup \{(i,j) \in \mathbb{Z}_{2(k+1)} \times \mathbb{Z}_{2(m+1)}: \text{\rm $i$ and $j$ are both odd}\}. 
\end{align*}
For each pair $(i_0,j_0) \in A$, we denote by $\rho_{i_0,j_0}$ the reflection across the geodesic arc $P_{i_0} Q_{j_0}$. Furthermore, we denote by $G$ the subgroup of $\text{\rm SO}(4)$ generated by the reflections $\{\rho_{i_0,j_0}: (i_0,j_0) \in A\}$. It is easy to see that each of the sets $\{P_i: \text{\rm $i$ is even}\}$, $\{P_i: \text{\rm $i$ is odd}\}$, $\{Q_j: \text{\rm $j$ is even}\}$, and $\{Q_j: \text{\rm $j$ is odd}\}$ is invariant under $G$. Hence, for each pair $(i,j) \in A_{\text{\rm even}}$, there exists a unique element $T_{i,j} \in G$ which maps the set $\{P_0,Q_0,P_1,Q_1\}$ to the set $\{P_i,Q_j,P_{i+1},Q_{j+1}\}$. Moreover, we have $\rho_{i_0,j_0} \circ T_{i,j} = T_{2i_0-i-1,2j_0-j-1}$ for all pairs $(i_0,j_0) \in A$ and $(i,j) \in A_{\text{\rm even}}$.

For each pair $(i,j) \in A_{\text{\rm even}}$, we denote by $\Gamma_{i,j}$ the geodesic quadrilateral with vertices $P_i$, $Q_j$, $P_{i+1}$, and $Q_{j+1}$. Moreover, we define 
\begin{align*} 
D_{i,j} 
&= \Big \{ x \in S^3: x_1 \, \sin \frac{\pi i}{k+1} < x_2 \, \cos \frac{\pi i}{k+1} \Big \} \\ 
&\cap \Big \{ x \in S^3: x_1 \, \sin \frac{\pi (i+1)}{k+1} > x_2 \, \cos \frac{\pi (i+1)}{k+1} \Big \} \\ 
&\cap \Big \{ x \in S^3: x_3 \, \sin \frac{\pi j}{m+1} < x_4 \, \cos \frac{\pi j}{m+1} \Big \} \\ 
&\cap \Big \{ x \in S^3: x_3 \, \sin \frac{\pi (j+1)}{m+1} > x_4 \, \cos \frac{\pi (j+1)}{m+1} \Big \}. 
\end{align*} 
The boundary of $D_{i,j}$ consists of four faces, each of which is totally geodesic. Thus, $D_{i,j}$ is a geodesic tetrahedron with vertices $P_i$, $Q_j$, $P_{i+1}$, and $Q_{j+1}$. In particular, we have $\Gamma_{i,j} \subset \partial D_{i,j}$ for each pair $(i,j) \in A_{\text{\rm even}}$. 

By Theorem 1 in \cite{Meeks-Yau}, there exists an embedded least area disk $\Sigma_{0,0} \subset \overline{D}_{0,0}$ with boundary $\partial \Sigma_{0,0} = \Gamma_{0,0}$. For each pair $(i,j) \in A_{\text{\rm even}}$, we denote by $\Sigma_{i,j} \subset \overline{D}_{i,j}$ the image of $\Sigma_{0,0}$ under the map $T_{i,j} \in G$. Clearly, $\Sigma_{i,j}$ is an embedded minimal disk in $\overline{D}_{i,j}$ with boundary $\partial \Sigma_{i,j} = \Gamma_{i,j}$. Moreover, since $\rho_{i_0,j_0} \circ T_{i,j} = T_{2i_0-i-1,2j_0-j-1}$, the reflection $\rho_{i_0,j_0}$ maps $\Sigma_{i,j}$ to $\Sigma_{2i_0-i-1,2j_0-j-1}$. Consequently, the union 
\[\Sigma = \bigcup_{(i,j) \in A_{\text{\rm even}}} \Sigma_{i,j}\] 
is invariant under the group $G$. Moreover, $\Sigma$ is a minimal surface away from the geodesic arcs $P_{i_0} Q_{j_0}$, where $(i_0,j_0) \in A$, and the density of $\Sigma$ along the geodesic arc $P_{i_0} Q_{i_0}$ is equal to $1$. Since $\Sigma$ is invariant under the reflection $\rho_{i_0,j_0}$, we conclude that $\Sigma$ is smooth away from the set $\{P_i: i \in \mathbb{Z}_{2(k+1)}\} \cup \{Q_j: j \in \mathbb{Z}_{2(m+1)}\}$. Using the removable singularities theorem for harmonic maps with finite energy, we conclude that $\Sigma$ is smooth. Finally, since each surface $\Sigma_{i,j}$ is embedded and the cells $D_{i,j}$ are disjoint, it follows that the surface $\Sigma$ is embedded as well. 

Finally, let us compute the genus of $\Sigma$. The geodesic quadrilateral $\Gamma_{i,j}$ has interior angles $\frac{\pi}{m+1}$, $\frac{\pi}{k+1}$, $\frac{\pi}{m+1}$, and $\frac{\pi}{k+1}$. Since $\Sigma_{i,j}$ is homeomorphic to a disk, the Gauss-Bonnet theorem implies that 
\[-\int_{\Sigma_{i,j}} K = 2\pi - \frac{\pi}{m+1} - \frac{\pi}{k+1} - \frac{\pi}{m+1} - \frac{\pi}{k+1} = \frac{2\pi (km-1)}{(k+1)(m+1)}.\] 
Since the set $A_{\text{\rm even}}$ has cardinality $2(k+1)(m+1)$, we conclude that 
\[4\pi(g-1) = -\int_\Sigma K = -\sum_{(i,j) \in A_{\text{\rm even}}} \int_{\Sigma_{i,j}} K = 4\pi (km-1),\] 
where $g$ denotes the genus of $\Sigma$. Thus, the surface $\Sigma$ has genus $g=km$. \\

\begin{theorem}[H.~Karcher, U.~Pinkall, and I.~Sterling \cite{Karcher-Pinkall-Sterling}]
There exist additional examples of embedded minimal surfaces in $S^3$, which are not part of the family obtained by Lawson. These surfaces have genus $3$, $5$, $6$, $7$, $11$, $19$, $73$, and $601$.
\end{theorem}

The construction in \cite{Karcher-Pinkall-Sterling} is similar in spirit to Lawson's construction; it uses tesselations of $S^3$ into cells that have the symmetry of a Platonic solid in $\mathbb{R}^3$. 

Very recently, Choe and Soret \cite{Choe-Soret2} announced a new construction of embedded minimal surfaces in $S^3$ which are obtained by desingularizing a union of Clifford tori. The proof of Choe and Soret is inspired by Lawson's construction, and uses reflection symmetries in a crucial way. There is an alternative construction by Kapouleas and Wiygul \cite{Kapouleas-Wiygul} which relies on gluing techniques and the implicit function theorem.

\begin{theorem}[J.~Choe, M.~Soret \cite{Choe-Soret2}]
\label{desingularization}
There exists a family of embedded minimal surfaces $\Sigma_{m,l}$ in $S^3$ of genus $1+4m(m-1)l$. Moreover, the surface $\Sigma_{m,l}$ can be viewed as a desingularization of the union $\bigcup_{j=0}^{m-1} T_j$, where 
\[T_j = \Big \{ x \in S^3: (x_1x_4+x_2x_3) \cos \frac{\pi j}{m} = (x_1x_3-x_2x_4) \sin \frac{\pi j}{m} \Big \}.\] 
\end{theorem} 

\textbf{Sketch of the proof of Theorem \ref{desingularization}.} 
Let us define 
\[P_i = \Big ( \cos \frac{\pi i}{2lm},\sin \frac{\pi i}{2lm},0,0 \Big )\] 
and 
\[Q_{i,j} = \frac{1}{\sqrt{2}} \, \Big ( \cos \frac{\pi i}{2lm},\sin \frac{\pi i}{2lm},\cos \frac{\pi (i - 2lj)}{2lm},-\sin \frac{\pi (i-2lj)}{2lm} \Big )\]
for $i \in \mathbb{Z}_{4ml}$ and $j \in \mathbb{Z}_{2m}$. 
Moreover, let 
\begin{align*} 
A_{\text{\rm even}} 
&= \{(i,j) \in \mathbb{Z}_{4ml} \times \mathbb{Z}_{2m}: \text{\rm $i$ and $j$ are both even}\} \\ 
&\cup \{(i,j) \in \mathbb{Z}_{4ml} \times \mathbb{Z}_{2m}: \text{\rm $i$ and $j$ are both odd}\}. 
\end{align*}
For each pair $(i,j) \in A_{\text{\rm even}}$, we denote by $\Delta_{i,j}$ the geodesic polygon with vertices $Q_{i,j} P_i Q_{i,j+1} Q_{i+1,j+1} P_{i+1} Q_{i+1,j}$. Note that the geodesic arc $P_i Q_{i,j}$ is contained in the intersection 
\[T_j \cap \Big \{ x \in S^3: x_1 \sin \frac{\pi i}{2ml} = x_2 \cos \frac{\pi i}{2ml} \Big \}.\] 
Moreover, the geodesic arc $Q_{i,j} Q_{i+1,j}$ is contained in the intersection 
\[T_j \cap \{x \in S^3: x_1^2+x_2^2=x_3^2+x_4^2\}.\] 
Given any pair $(i,j) \in A_{\text{\rm even}}$, we define 
\begin{align*} 
U_{i,j} &= \Big \{ x \in S^3: x_1 \sin \frac{\pi i}{2ml} < x_2 \cos \frac{\pi i}{2ml} \Big \} \\ 
&\cap \Big \{ x \in S^3: x_1 \sin \frac{\pi (i+1)}{2ml} > x_2 \cos \frac{\pi (i+1)}{2ml} \Big \} \\ 
&\cap \Big \{ x \in S^3: (x_1x_4+x_2x_3) \cos \frac{\pi j}{m} > (x_1x_3-x_2x_4) \sin \frac{\pi j}{m} \Big \} \\ 
&\cap \Big \{ x \in S^3: (x_1x_4+x_2x_3) \cos \frac{\pi (j+1)}{m} < (x_1x_3-x_2x_4) \sin \frac{\pi (j+1)}{m} \Big \} \\ 
&\cap \{ x \in S^3: x_1^2+x_2^2 > x_3^2+x_4^2\}. 
\end{align*}
Note that $U_{i,j}$ is a mean convex domain in $S^3$. In fact, the boundary of $U_{i,j}$ consists of five faces: two of these faces are totally geodesic, and each of the other three faces is congruent to a piece of the Clifford torus. Moreover, it is straightforward to verify that $\Delta_{i,j} \subset \partial U_{i,j}$.

Let $S_{0,0} \subset \overline{U}_{0,0}$ be an embedded least area disk with boundary $\partial S_{0,0} = \Delta_{0,0}$. The reflection across the geodesic arc $P_{i_0} Q_{i_0,j_0}$ maps the region $U_{i,j}$ to $U_{2i_0-i-1,2j_0-j-1}$ and the polygon $\Delta_{i,j}$ to $\Delta_{2i_0-i-1,2j_0-j-1}$. Hence, by successive reflection across geodesic arcs on the boundary, one obtains a family of embedded least area disks $S_{i,j} \subset \overline{U}_{i,j}$ with boundary $\partial S_{i,j} = \Delta_{i,j}$. The union 
\[S = \bigcup_{(i,j) \in A_{\text{\rm even}}} S_{i,j}\] 
is a smooth minimal surface which is contained in the region $\{x \in S^3: x_1^2+x_2^2 \geq x_3^2+x_4^2\}$. Moreover, the boundary of $S$ lies on the Clifford torus $\{x \in S^3: x_1^2+x_2^2=x_3^2+x_4^2\}$. Choe and Soret then show that the union 
\[\Sigma = S \cup \{(x_3,-x_4,x_1,-x_2): (x_1,x_2,x_3,x_4) \in S\}\] 
is a smooth minimal surface in $S^3$. This surface is clearly embedded. 

It remains to compute the genus of $\Sigma$. The interior angles of the geodesic polygon $\Delta_{i,j}$ are $\frac{\pi}{2}$, $\frac{\pi}{m}$, $\frac{\pi}{2}$, $\frac{\pi}{2}$, $\frac{\pi}{m}$, and $\frac{\pi}{2}$. Therefore, 
\[-\int_{S_{i,j}} K = 4\pi - \frac{\pi}{2} - \frac{\pi}{m} - \frac{\pi}{2} - \frac{\pi}{2} - \frac{\pi}{m} - \frac{\pi}{2} = \frac{2\pi (m-1)}{m}\] 
by the Gauss-Bonnet theorem. Since the set $A_{\text{\rm even}}$ has cardinality $4m^2l$, it follows that 
\[4\pi(g-1) = -\int_\Sigma K = -2 \sum_{(i,j) \in A_{\text{\rm even}}} \int_{S_{i,j}} K = 16\pi m(m-1)l.\] 
Consequently, $g = 1+4m(m-1)l$, as claimed. \\

In the remainder of this section, we describe another family of embedded minimal surfaces in $S^3$, which was constructed by Kapouleas and Yang \cite{Kapouleas-Yang} using gluing techniques. The idea here is to take two nearby copies of the Clifford torus, and join them by a large number of catenoid bridges. In this way, one obtains a family of approximate solutions of the minimal surface equation, and Kapouleas and Yang showed that these surfaces can be deformed to exact solutions of the minimal surface equation by means of the implicit function theorem. 

In the following, we sketch the construction of the initial surfaces in \cite{Kapouleas-Yang}. The Clifford torus can be parametrized by a map $F: \mathbb{R}^2 \to S^3$, where 
\[F(s,t) = \frac{1}{\sqrt{2}} \, (\cos(\sqrt{2} \, s),\sin(\sqrt{2} \, s),\cos(\sqrt{2} \, t),\sin(\sqrt{2} \, t)).\] 
The map $F$ can be extended to a map $\Phi: \mathbb{R}^2 \times (-\frac{\pi}{4},\frac{\pi}{4}) \to S^3$ by 
\begin{align*} 
\Phi(s,t,u) 
&= \sin \Big ( u+\frac{\pi}{4} \Big ) \, (\cos(\sqrt{2} \, s),\sin(\sqrt{2} \, s),0,0) \\ 
&+ \cos \Big ( u+\frac{\pi}{4} \Big ) \, (0,0,\cos(\sqrt{2} \, t),\sin(\sqrt{2} \, t)) 
\end{align*}
(see \cite{Kapouleas-Yang}, equation (2.1)). Note that $F(s,t) = \Phi(s,t,0)$. Moreover, the pull-back of the round metric on $S^3$ under the map $\Phi$ can be expressed as 
\[g = (1 + \sin(2u)) \, ds \otimes ds + (1 - \sin(2u)) \, dt \otimes dt + du \otimes du.\] 
The approximate solutions constructed in \cite{Kapouleas-Yang} depend on two parameters, an integer $m$ (which is assumed to be very large) and a real number $\zeta$ (which lies in a bounded interval). Following \cite{Kapouleas-Yang}, we put 
\[\tau = \frac{1}{m} \, e^{-\frac{m^2}{4\pi}+\zeta}.\] 
Let $\psi: \mathbb{R} \to [0,1]$ be a smooth cutoff function such that $\psi = 1$ on $(-\infty,1]$ and $\psi = 0$ on $[2,\infty)$. Kapouleas and Yang then define 
\[M_{\text{\rm cat}} = \bigg \{ \Phi(s,t,u): \text{\rm $\tau \leq \sqrt{s^2+t^2} \leq \frac{1}{m}$ and $\frac{|u|}{\tau} = \text{\rm arcosh} \, \frac{\sqrt{s^2+t^2}}{\tau}$} \bigg \}\] 
and 
\begin{align*} 
M_{\text{\rm tor}} &= \bigg \{ \Phi(s,t,u): \text{\rm $\sqrt{s^2+t^2} \geq \frac{1}{m}$, $\max\{|s|,|t|\} \leq \frac{\pi}{\sqrt{2} \, m}$,} \\ 
&\hspace{25mm} \text{\rm and $\frac{|u|}{\tau} = \psi(m\sqrt{s^2+t^2}) \, \text{\rm arcosh} \, \frac{\sqrt{s^2+t^2}}{\tau}$} \\ 
&\hspace{35mm} + (1-\psi(m\sqrt{s^2+t^2})) \, \text{\rm arcosh} \, \frac{1}{m\tau} \bigg \}. 
\end{align*} 
The union $M = M_{\text{\rm cat}} \cup M_{\text{\rm tor}}$ is a smooth surface with boundary. By gluing together $m^2$ rotated copies of the surface $M$, we obtain a closed, embedded surface in $S^3$ of genus $m^2+1$. This surface depends on the parameters $m$ and $\zeta$, and will be denoted by $\Sigma_{m,\zeta}$. Since the catenoid in $\mathbb{R}^3$ has zero mean curvature, the mean curvature of the surface $\Sigma_{m,\zeta}$ is small when $m$ is sufficiently large (see \cite{Kapouleas-Yang}, Lemma 3.18, for a precise statement). 

The key issue is to deform the surface $\Sigma_{m,\zeta}$ to an exact solution of the minimal surface equation. This is a difficult problem, since the linearized operator has a non-trivial kernel. Taking into account the symmetries of the problem, the approximate kernel turns out to be one-dimensional. In fact, Kapouleas and Yang show that the approximate kernel stems from the constant functions on $M_{\text{\rm tor}}$ (see \cite{Kapouleas-Yang}, Proposition 4.14), and this obstacle can be overcome by a suitable choice of the parameter $\zeta$: 

\begin{theorem}[N.~Kapouleas, S.D.~Yang \cite{Kapouleas-Yang}]
\label{gluing}
If $m$ is sufficiently large, then there exists a real number $\zeta_m$ with the property that $\Sigma_{m,\zeta_m}$ can be deformed to an embedded minimal surface $\hat{\Sigma}_m$ of genus $m^2+1$. 
\end{theorem}

While the construction of Kapouleas and Yang is not explicit, the estimates in \cite{Kapouleas-Yang} provide a very precise description of the surfaces $\hat{\Sigma}_m$ when $m$ is large. In particular, the surfaces $\hat{\Sigma}_m$ converge, in the sense of varifolds, to the Clifford torus with multiplicity $2$ as $m \to \infty$. Finally, we note that Kapouleas \cite{Kapouleas-survey} has recently announced a similar doubling construction for the equator.

\section{Uniqueness questions for minimal surfaces and the Lawson conjecture}

In this section, we discuss uniqueness results for minimal surfaces of genus $0$ and $1$. In 1966, Almgren proved the following uniqueness theorem in the genus $0$ case: 

\begin{theorem}[F.J.~Almgren, Jr. \cite{Almgren}]
\label{almgren.genus.0}
The equator is the only immersed minimal surface in $S^3$ of genus $0$ (up to rigid motions in $S^3$).
\end{theorem}

\textbf{Proof of Theorem \ref{almgren.genus.0}.} The proof relies on a Hopf differential argument. To explain this, let $F: S^2 \to S^3$ be a conformal minimal immersion, and let $h$ denote its second fundamental form. We will identify $S^2$ with $\mathbb{C} \cup \{\infty\}$, where the north pole on $S^2$ corresponds to the point at infinity. It follows from the Codazzi equations that the function $h \big ( \frac{\partial}{\partial z},\frac{\partial}{\partial z} \big )$ is holomorphic. Moreover, since the immersion $F$ is smooth at the north pole, the function $h \big ( \frac{\partial}{\partial z},\frac{\partial}{\partial z} \big )$ vanishes at the north pole. By Liouville's theorem, the function $h \big ( \frac{\partial}{\partial z},\frac{\partial}{\partial z} \big )$ vanishes identically. On the other hand, we have $h \big ( \frac{\partial}{\partial z},\frac{\partial}{\partial \bar{z}} \big ) = 0$ since the mean curvature of $F$ vanishes. Thus, $F$ is totally geodesic, hence congruent to the equator. \\

In 1970, Lawson conjectured a similar uniqueness property for minimal tori in $S^3$. Specifically, Lawson conjectured the following: 

\begin{conjecture}[H.B.~Lawson, Jr. \cite{Lawson3}]
The Clifford torus is the only embedded minimal surface in $S^3$ of genus $1$ (up to rigid motions in $S^3$).
\end{conjecture} 

Note that Lawson's conjecture is false if the surface if we allow the surface to have self-intersections (see \cite{Lawson1} or Theorem \ref{alexandrov.immersed.tori} above).

In March 2012, we gave an affirmative answer to Lawson's conjecture (cf. \cite{Brendle2}). One of the main difficulties is that any proof of Lawson's conjecture has to exploit the assumption that $\Sigma$ is embedded, as well as the condition that $\Sigma$ has genus $1$. In order to exploit the latter condition, we make use of the following result due to Lawson:

\begin{proposition}[H.B.~Lawson, Jr. \cite{Lawson2}]
\label{absence.of.umbilic.points}
An immersed minimal surface in $S^3$ of genus $1$ has no umbilic points; in other words, the second fundamental form is non-zero at each point on the surface.
\end{proposition}

\textbf{Proof of Proposition \ref{absence.of.umbilic.points}.} Let $F: \Sigma \to S^3$ be a conformal minimal immersion of genus $1$, and let $h$ denote its second fundamental form. We may write $\Sigma = \mathbb{C} / \Lambda$, where $\Lambda$ is a lattice in $\mathbb{C}$. As above, the Codazzi equations imply that the expression $h \big ( \frac{\partial}{\partial z},\frac{\partial}{\partial z} \big )$ defines a holomorphic function on $\mathbb{C} / \Lambda$. By Liouville's theorem, we have $h \big ( \frac{\partial}{\partial z},\frac{\partial}{\partial z} \big ) = c$ for some constant $c$. If $c = 0$, then the surface is a totally geodesic two-sphere, contradicting our assumption that $\Sigma$ has genus $1$. Thus, $c \neq 0$, and the second fundamental form is non-zero at each point on the surface. This completes the proof of Proposition \ref{absence.of.umbilic.points}. \\

Moreover, we will need the following result, which is a consequence of the well-known Simons identity (cf. \cite{Simons}):

\begin{proposition} 
\label{simons.identity}
Suppose that $F: \Sigma \to S^3$ is an embedded minimal torus in $S^3$. Then the norm of the second fundamental form satisfies the partial differential equation 
\[\Delta_\Sigma(|A|) - \frac{\big | \nabla |A| \big |^2}{|A|} + (|A|^2 - 2) \, |A| = 0.\] 
\end{proposition}

\textbf{Sketch of the proof of Proposition \ref{simons.identity}.} 
Using the Simons identity 
\[\Delta_\Sigma(|A|^2) - 2 \, |\nabla A|^2 + 2 \, (|A|^2-2) \, |A|^2 = 0,\] 
we obtain 
\[\Delta_\Sigma(|A|) + \frac{\big | \nabla |A| \big |^2}{|A|} - \frac{|\nabla A|^2}{|A|} + (|A|^2-2) \, |A| = 0.\] 
On the other hand, the Codazzi equations imply that $|\nabla A|^2 = 2 \, \big | \nabla |A| \big |^2$. From this, the assertion follows. \\

The proof of the Lawson conjecture in \cite{Brendle2} involves an application of the maximum principle to a function that depends on a pair of points. This technique was pioneered by Huisken \cite{Huisken} in his work on the curve shortening flow for embedded curves in the plane. Specifically, Huisken was able to give a lower bound for the chord distance in terms of the arc length. This gives a new proof of Grayson's theorem, which asserts that any embedded curve shrinks to a point in finite time and becomes round after rescaling (cf. \cite{Grayson}, \cite{Hamilton}). Using a similar method, Andrews \cite{Andrews} obtained an alternative proof of the noncollapsing property for mean curvature flow. The noncollapsing theorem for the mean curvature flow was first stated in a paper by Sheng and Wang in \cite{Sheng-Wang}; the result is a direct consequence of the work of White on the structure of singularities in the mean curvature flow (cf. \cite{White1}, \cite{White2}, \cite{White3}).

The argument in \cite{Brendle2} uses a different quantity, which involves the norm of the second fundamental form. A major difficulty we encounter in this approach is that the Simons identity for the norm of the second fundamental form contains a gradient term, which turns out to have an unfavorable sign. As a result, the calculation becomes extremely subtle and we need to make use of every available piece of information. We will describe the details below. Suppose that $F: \Sigma \to S^3$ is a minimal immersion of a genus $1$ surface into $S^3$. Moreover, let $\nu(x) \in T_{F(x)} S^3$ be a unit normal vector field. For abbreviation, we define a smooth function $\Psi: \Sigma \to \mathbb{R}$ by 
\[\Psi(x) = \frac{1}{\sqrt{2}} \, |A(x)|,\] 
where $|A(x)|$ denotes the norm of the second fundamental form. Since $F$ is a minimal immersion, the principal curvatures at the point $x$ satisfy $|\lambda_1| = |\lambda_2| = \Psi(x)$. Note that the function $\Psi$ is strictly positive by Proposition \ref{absence.of.umbilic.points}. 

Given any number $\alpha \geq 1$, we define a function $Z_\alpha: \Sigma \times \Sigma \to \mathbb{R}$ by 
\begin{equation} 
Z_\alpha(x,y) = \alpha \, \Psi(x) \, (1 - \langle F(x),F(y) \rangle) + \langle \nu(x),F(y) \rangle. 
\end{equation}
We begin by compute the gradient of the function $Z_\alpha$. To that end, we fix two distinct points $\bar{x},\bar{y} \in \Sigma$. Moreover, let $(x_1,x_2)$ be a system of geodesic normal coordinates around $\bar{x}$, and let $(y_1,y_2)$ be a geodesic normal coordinates around $\bar{y}$. Without loss of generality, we may assume that the second fundamental form at $\bar{x}$ is diagonal, so that $h_{11}(\bar{x}) = \lambda_1$, $h_{12}(\bar{x}) = 0$, and $h_{22}(\bar{x}) = \lambda_2$. 

The first derivatives of the function $Z_\alpha$ are given by 
\begin{align} 
\label{gradient.a}
\frac{\partial Z_\alpha}{\partial x_i}(\bar{x},\bar{y}) 
&= \alpha \, \frac{\partial \Psi}{\partial x_i}(\bar{x}) \, (1 - \langle F(\bar{x}),F(\bar{y}) \rangle) \notag \\ 
&- \alpha \, \Psi(\bar{x}) \, \Big \langle \frac{\partial F}{\partial x_i}(\bar{x}),F(\bar{y}) \Big \rangle + h_i^k(\bar{x}) \, \Big \langle \frac{\partial F}{\partial x_k}(\bar{x}),F(\bar{y}) \Big \rangle 
\end{align}
and 
\begin{equation} 
\label{gradient.b}
\frac{\partial Z_\alpha}{\partial y_i}(\bar{x},\bar{y}) = -\alpha \, \Psi(\bar{x}) \, \Big \langle F(\bar{x}),\frac{\partial F}{\partial y_i}(\bar{y}) \Big \rangle + \Big \langle \nu(\bar{x}),\frac{\partial F}{\partial y_i}(\bar{y}) \Big \rangle. 
\end{equation}
We next consider the second order derivatives of $Z$ at the point $(\bar{x},\bar{y})$. 

\begin{lemma}
\label{laplacian}
The Laplacian of $Z_\alpha$ with respect to $x$ satisfies an inequality of the form 
\begin{align} 
\label{second.order.condition.a}
&\sum_{i=1}^2 \frac{\partial^2 Z_\alpha}{\partial x_i^2}(\bar{x},\bar{y}) \notag \\ 
&\leq 2\alpha \, \Psi(\bar{x}) - \frac{\alpha^2-1}{\alpha} \, \frac{\Psi(\bar{x})}{1-\langle F(\bar{x}),F(\bar{y}) \rangle} \sum_{i=1}^2 \Big \langle \frac{\partial F}{\partial x_i}(\bar{x}),F(\bar{y}) \Big \rangle^2 \\ 
&+ \Lambda_1(|F(\bar{x})-F(\bar{y})|) \, \bigg ( |Z_\alpha(\bar{x},\bar{y})| + \sum_{i=1}^2 \Big | \frac{\partial Z_\alpha}{\partial x_i}(\bar{x},\bar{y}) \Big | \bigg ), \notag
\end{align} 
where $\Lambda_1: (0,\infty) \to (0,\infty)$ is a continuous function. Moreover, the Laplacian of $Z_\alpha$ with respect to $y$ satisfies 
\begin{equation} 
\label{second.order.condition.b}
\sum_{i=1}^2 \frac{\partial^2 Z_\alpha}{\partial y_i^2}(\bar{x},\bar{y}) \leq 2\alpha \, \Psi(\bar{x}) + 2 \, |Z_\alpha(\bar{x},\bar{y})|. 
\end{equation}
\end{lemma}

\textbf{Proof of Lemma \ref{laplacian}.} 
By the Codazzi equations, we have 
\[\sum_{i=1}^2 \frac{\partial}{\partial x_i} h_i^k(\bar{x}) = 0.\] 
Using this identity, we compute 
\begin{align*} 
&\sum_{i=1}^2 \frac{\partial^2 Z_\alpha}{\partial x_i^2}(\bar{x},\bar{y}) \\ 
&= \alpha \sum_{i=1}^2 \frac{\partial^2 \Psi}{\partial x_i^2}(\bar{x}) \, (1 - \langle F(\bar{x}),F(\bar{y}) \rangle) - 2\alpha \sum_{i=1}^2 \frac{\partial \Psi}{\partial x_i}(\bar{x}) \, \Big \langle \frac{\partial F}{\partial x_i}(\bar{x}),F(\bar{y}) \Big \rangle \\ 
&+ 2\alpha \, \Psi(\bar{x}) \, \langle F(\bar{x}),F(\bar{y}) \rangle - |A(\bar{x})|^2 \, \langle \nu(\bar{x}),F(\bar{y}) \rangle \\ 
&= \alpha \, \big ( \Delta_\Sigma \Psi(\bar{x}) + (|A(\bar{x})|^2-2) \, \Psi(\bar{x}) \big ) \, (1 - \langle F(\bar{x}),F(\bar{y}) \rangle) + 2\alpha \, \Psi(\bar{x}) \\ 
&- 2\alpha \sum_{i=1}^2 \frac{\partial \Psi}{\partial x_i}(\bar{x}) \, \Big \langle \frac{\partial F}{\partial x_i}(\bar{x}),F(\bar{y}) \Big \rangle - |A(\bar{x})|^2 \, Z_\alpha(\bar{x},\bar{y}). 
\end{align*} 
Proposition \ref{simons.identity} implies that 
\[\Delta_\Sigma \Psi - \frac{|\nabla \Psi|^2}{\Psi} + (|A|^2 - 2) \, \Psi = 0.\] 
This gives 
\begin{align*} 
&\sum_{i=1}^2 \frac{\partial^2 Z_\alpha}{\partial x_i^2}(\bar{x},\bar{y}) \\ 
&= \alpha \, \frac{|\nabla \Psi(\bar{x})|^2}{\Psi(\bar{x})} \, (1 - \langle F(\bar{x}),F(\bar{y}) \rangle) + 2\alpha \, \Psi(\bar{x}) \\ 
&- 2\alpha \sum_{i=1}^2 \frac{\partial \Psi}{\partial x_i}(\bar{x}) \, \Big \langle \frac{\partial F}{\partial x_i}(\bar{x}),F(\bar{y}) \Big \rangle - |A(\bar{x})|^2 \, Z_\alpha(\bar{x},\bar{y}). 
\end{align*} 
The expression on the right hand side can be rewritten as 
\begin{align*} 
&\sum_{i=1}^2 \frac{\partial^2 Z_\alpha}{\partial x_i^2}(\bar{x},\bar{y}) \\ 
&= \frac{\alpha}{\Psi(\bar{x}) \, (1-\langle F(\bar{x}),F(\bar{y}) \rangle)} \\ 
&\hspace{15mm} \cdot \sum_{i=1}^2 \bigg ( \frac{\partial \Psi}{\partial x_i}(\bar{x}) \, (1 - \langle F(\bar{x}),F(\bar{y}) \rangle) - \Psi(\bar{x}) \, \Big \langle \frac{\partial F}{\partial x_i}(\bar{x}),F(\bar{y}) \Big \rangle \bigg )^2 \\ 
&+ 2\alpha \, \Psi(\bar{x}) - \frac{\alpha \, \Psi(\bar{x})}{1-\langle F(\bar{x}),F(\bar{y}) \rangle} \sum_{i=1}^2 \Big \langle \frac{\partial F}{\partial x_i}(\bar{x}),F(\bar{y}) \Big \rangle^2 - |A(\bar{x})|^2 \, Z_\alpha(\bar{x},\bar{y}). 
\end{align*} 
Using the relation (\ref{gradient.a}), we conclude that 
\begin{align*} 
\sum_{i=1}^2 \frac{\partial^2 Z_\alpha}{\partial x_i^2}(\bar{x},\bar{y}) 
&\leq \frac{1}{\alpha \, \Psi(\bar{x}) \, (1-\langle F(\bar{x}),F(\bar{y}) \rangle)} \sum_{i=1}^2 \lambda_i^2 \, \Big \langle \frac{\partial F}{\partial x_i}(\bar{x}),F(\bar{y}) \Big \rangle^2 \\ 
&+ 2\alpha \, \Psi(\bar{x}) - \frac{\alpha \, \Psi(\bar{x})}{1-\langle F(\bar{x}),F(\bar{y}) \rangle} \sum_{i=1}^2 \Big \langle \frac{\partial F}{\partial x_i}(\bar{x}),F(\bar{y}) \Big \rangle^2 \\ 
&+ \Lambda_1(|F(\bar{x}) - F(\bar{y})|) \, \bigg ( |Z_\alpha(\bar{x},\bar{y})| + \sum_{i=1}^2 \Big | \frac{\partial Z_\alpha}{\partial x_i}(\bar{x},\bar{y}) \Big | \bigg ), 
\end{align*} 
where $\Lambda_1: (0,\infty) \to (0,\infty)$ is a continuous function. Since $\lambda_1^2 = \lambda_2^2 = \Psi(\bar{x})^2$, the identity (\ref{second.order.condition.a}) follows. Finally, we have 
\begin{align*} 
\sum_{i=1}^2 \frac{\partial^2 Z_\alpha}{\partial y_i^2}(\bar{x},\bar{y}) 
&= 2\alpha \, \Psi(\bar{x}) \, \langle F(\bar{x}),F(\bar{y}) \rangle - 2 \, \langle \nu(\bar{x}),F(\bar{y}) \rangle \\ 
&= 2\alpha \, \Psi(\bar{x}) - 2 \, Z_\alpha(\bar{x},\bar{y}). 
\end{align*}
This proves (\ref{second.order.condition.b}). \\

Finally, we estimate the mixed partial derivatives of $Z_\alpha$.

\begin{lemma}
\label{mixed.partials}
For a suitable choice of the coordinate system $(y_1,y_2)$, we have 
\begin{align*} 
&\sum_{i=1}^2 \frac{\partial^2 Z_\alpha}{\partial x_i \, \partial y_i}(\bar{x},\bar{y}) \\ 
&\leq -2\alpha \, \Psi(\bar{x}) + \Lambda_4(|F(\bar{x})-F(\bar{y})|) \, \bigg ( |Z_\alpha(\bar{x},\bar{y})| + \sum_{i=1}^2 \Big | \frac{\partial Z_\alpha}{\partial x_i}(\bar{x},\bar{y}) \Big | + \sum_{i=1}^2 \Big | \frac{\partial Z_\alpha}{\partial y_i}(\bar{x},\bar{y}) \Big | \bigg ), 
\end{align*}
where $\Lambda_4: (0,\infty) \to (0,\infty)$ is a continuous function.
\end{lemma}

\textbf{Proof of Lemma \ref{mixed.partials}.} 
Let $w_i$ denote the reflection of the vector $\frac{\partial F}{\partial x_i}(\bar{x})$ across the hyperplane orthogonal to $F(\bar{x})-F(\bar{y})$, so that 
\[w_i = \frac{\partial F}{\partial x_i}(\bar{x}) - 2 \, \Big \langle \frac{\partial F}{\partial x_i}(\bar{x}),\frac{F(\bar{x})-F(\bar{y})}{|F(\bar{x})-F(\bar{y})|} \Big \rangle \, \frac{F(\bar{x})-F(\bar{y})}{|F(\bar{x})-F(\bar{y})|}.\] 
If $Z_\alpha(\bar{x},\bar{y}) = 0$ and $\frac{\partial Z_\alpha}{\partial y_i}(\bar{x},\bar{y}) = 0$, then we have 
\[\text{\rm span} \Big \{ \frac{\partial F}{\partial y_1}(\bar{y}),\frac{\partial F}{\partial y_2}(\bar{y}) \Big \} = \text{\rm span}\{w_1,w_2\}.\] 
Hence, in this case, we may choose the coordinate system $(y_1,y_2)$ so that $\frac{\partial F}{\partial y_i}(\bar{y}) = w_i$ for $i=1,2$. 

We now return to the general case. We may choose the coordinate system $(y_1,y_2)$ in such a way that 
\[\sum_{i=1}^2 \Big | \frac{\partial F}{\partial y_i}(\bar{y}) - w_i \Big | \leq \Lambda_2(|F(\bar{x})-F(\bar{y})|) \, \bigg ( |Z_\alpha(\bar{x},\bar{y})| + \sum_{i=1}^2 \Big | \frac{\partial Z_\alpha}{\partial y_i}(\bar{x},\bar{y}) \Big | \bigg ),\] 
where $\Lambda_2: (0,\infty) \to (0,\infty)$ is a continuous function. For this choice of the coordinate system $(y_1,y_2)$, we have 
\begin{align*} 
&\frac{\partial^2 Z_\alpha}{\partial x_i \, \partial y_i}(\bar{x},\bar{y}) \\ 
&= -\alpha \, \frac{\partial \Psi}{\partial x_i}(\bar{x}) \, \Big \langle F(\bar{x}),\frac{\partial F}{\partial y_i}(\bar{y}) \Big \rangle + (\lambda_i - \alpha \, \Psi(\bar{x})) \, \Big \langle \frac{\partial F}{\partial x_i}(\bar{x}),\frac{\partial F}{\partial y_i}(\bar{y}) \Big \rangle \\ 
&= (\lambda_i - \alpha \, \Psi(\bar{x})) \, \Big \langle \frac{\partial F}{\partial x_i}(\bar{x}),\frac{\partial F}{\partial y_i}(\bar{y}) \Big \rangle \\ 
&+ \frac{1}{1 - \langle F(\bar{x}),F(\bar{y}) \rangle} \, (\lambda_i - \alpha \, \Psi(\bar{x})) \, \Big \langle \frac{\partial F}{\partial x_i}(\bar{x}),F(\bar{y}) \Big \rangle \, \Big \langle F(\bar{x}),\frac{\partial F}{\partial y_i}(\bar{y}) \Big \rangle \\ 
&- \frac{1}{1 - \langle F(\bar{x}),F(\bar{y}) \rangle} \, \Big \langle F(\bar{x}),\frac{\partial F}{\partial y_i}(\bar{y}) \Big \rangle \, \frac{\partial Z_\alpha}{\partial x_i}(\bar{x},\bar{y}) \\ 
&= (\lambda_i - \alpha \, \Psi(\bar{x})) \, \Big \langle w_i,\frac{\partial F}{\partial y_i}(\bar{y}) \Big \rangle \\ 
&- \frac{1}{1 - \langle F(\bar{x}),F(\bar{y}) \rangle} \, \Big \langle F(\bar{x}),\frac{\partial F}{\partial y_i}(\bar{y}) \Big \rangle \, \frac{\partial Z_\alpha}{\partial x_i}(\bar{x},\bar{y}). 
\end{align*} 
Thus, we conclude that 
\begin{align*} 
&\frac{\partial^2 Z_\alpha}{\partial x_i \, \partial y_i}(\bar{x},\bar{y}) \\ 
&\leq \lambda_i - \alpha \, \Psi(\bar{x}) \\ 
&+ \Lambda_3(|F(\bar{x})-F(\bar{y})|) \, \bigg ( |Z_\alpha(\bar{x},\bar{y})| + \sum_{i=1}^2 \Big | \frac{\partial Z_\alpha}{\partial x_i}(\bar{x},\bar{y}) \Big | + \sum_{i=1}^2 \Big | \frac{\partial Z_\alpha}{\partial y_i}(\bar{x},\bar{y}) \Big | \bigg ), 
\end{align*} 
where $\Lambda_3: (0,\infty) \to (0,\infty)$ is a continuous function. Hence, the assertion follows by summation over $i$. This completes the proof of Lemma \ref{mixed.partials}. \\

Combining Lemma \ref{laplacian} and Lemma \ref{mixed.partials}, we can draw the following conclusion: 

\begin{proposition}[S.~Brendle \cite{Brendle2}]
\label{key.estimate}
For a suitable choice of the coordinate system $(y_1,y_2)$, we have 
\begin{align*} 
&\sum_{i=1}^2 \frac{\partial^2 Z_\alpha}{\partial x_i^2}(\bar{x},\bar{y}) + 2 \sum_{i=1}^2 \frac{\partial^2 Z_\alpha}{\partial x_i \, \partial y_i}(\bar{x},\bar{y}) + \sum_{i=1}^2 \frac{\partial^2 Z_\alpha}{\partial y_i^2}(\bar{x},\bar{y}) \\ 
&\leq -\frac{\alpha^2-1}{\alpha} \, \frac{\Psi(\bar{x})}{1-\langle F(\bar{x}),F(\bar{y}) \rangle} \sum_{i=1}^2 \Big \langle \frac{\partial F}{\partial x_i}(\bar{x}),F(\bar{y}) \Big \rangle^2 \\ 
&+ \Lambda_5(|F(\bar{x})-F(\bar{y})|) \, \bigg ( |Z_\alpha(\bar{x},\bar{y})| + \sum_{i=1}^2 \Big | \frac{\partial Z_\alpha}{\partial x_i}(\bar{x},\bar{y}) \Big | + \sum_{i=1}^2 \Big | \frac{\partial Z_\alpha}{\partial y_i}(\bar{x},\bar{y}) \Big | \bigg ), 
\end{align*} 
where $\Lambda_5: (0,\infty) \to (0,\infty)$ is a continuous function.
\end{proposition}

After these preparations, we now state the main result in \cite{Brendle2}:

\begin{theorem}[S.~Brendle \cite{Brendle2}]
\label{lawson.conj} 
Let $F: \Sigma \to S^3$ be an embedded minimal surface in $S^3$ of genus $1$. Then $F$ is congruent to the Clifford torus.
\end{theorem}

\textbf{Proof of Theorem \ref{lawson.conj}.} 
Since $\Sigma$ is embedded and has no umbilic points, we have 
\[\kappa := \sup_{x,y \in \Sigma, \, x \neq y} \frac{|\langle \nu(x),F(y) \rangle|}{\Psi(x) \, (1-\langle F(x),F(y) \rangle} < \infty.\] 
We now distinguish two cases: 

\textit{Case 1:} Suppose first that $\kappa=1$. In this case, we have $Z_1(x,y) \geq 0$ for all points $x,y \in \Sigma$. Let us fix a point $\bar{x} \in \Sigma$, and let $\{e_1,e_2\}$ be an orthonormal basis of $T_{\bar{x}} \Sigma$ such that $h(e_1,e_1) = \Psi(\bar{x})$, $h(e_1,e_2) = 0$, and $h(e_2,e_2) = -\Psi(\bar{x})$. Moreover, we define 
\[\xi = \Psi(\bar{x}) \, F(\bar{x}) - \nu(\bar{x}) \in \mathbb{R}^4.\] 
Finally, we assume that $\sigma: \mathbb{R} \to \Sigma$ is a geodesic such that $\sigma(0) = \bar{x}$ and $\sigma'(0) = e_1$. The function 
\[f(t) = Z_1(\sigma(0),\sigma(t)) = \Psi(\bar{x}) - \langle \xi,F(\sigma(t)) \rangle\] 
is nonnegative for all $t$. A straightforward calculation gives 
\[f'(t) = -\langle \xi,dF_{\sigma(t)}(\sigma'(t)) \rangle,\]
\[f''(t) = \langle \xi,F(\sigma(t)) \rangle + h(\sigma'(t),\sigma'(t)) \, \langle \xi,\nu(\sigma(t)) \rangle,\] 
and 
\begin{align*} 
f'''(t) 
&= \langle \xi,dF_{\sigma(t)}(\sigma'(t)) \rangle + h(\sigma'(t),\sigma'(t)) \, \langle \xi,D_{\sigma'(t)} \nu \rangle \\ 
&+ (D_{\sigma'(t)}^\Sigma h)(\sigma'(t),\sigma'(t)) \, \langle \xi,\nu(\sigma(t)) \rangle. 
\end{align*} 
In particular, for $t=0$, we have $f(0) = f'(0) = f''(0) = 0$. Since the function $f(t)$ is nonnegative, we conclude that $f'''(0) = 0$. This implies that $(D_{e_1}^\Sigma h)(e_1,e_1) = 0$. Thus, $\langle \nabla \Psi(\bar{x}),e_1 \rangle = 0$. Replacing $\nu$ by $-\nu$, we obtain $\langle \nabla \Psi(\bar{x}),e_2 \rangle = 0$. Since the point $\bar{x}$ is arbitrary, the function $\Psi$ is constant, and the intrinsic Gaussian curvature of $\Sigma$ vanishes identically. By a result of Lawson \cite{Lawson1}, the surface is congruent to the Clifford torus. 

\textit{Case 2:} Suppose next that $\kappa > 1$. After replacing $\nu$ by $-\nu$, we may assume that \[\kappa = \sup_{x,y \in \Sigma, \, x \neq y} \Big ( -\frac{\langle \nu(x),F(y) \rangle}{\Psi(x) \, (1 - \langle F(x),F(y) \rangle)} \Big ).\] 
By definition of $\kappa$, the function $Z_\kappa$ is nonnegative, and the set 
\[\Omega = \{\bar{x} \in \Sigma: \text{\rm there exists a point $\bar{y} \in \Sigma \setminus \{\bar{x}\}$ such that $Z_\kappa(\bar{x},\bar{y}) = 0$}\}\] 
is non-empty. Using Proposition \ref{key.estimate} and Bony's maximum principle for degenerate elliptic equations (cf. \cite{Bony}), we conclude that the set $\Omega$ is open. 

We claim that $\nabla \Psi(\bar{x}) = 0$ for each point $\bar{x} \in \Omega$. Indeed, if $\bar{x} \in \Omega$, then we can find a point $\bar{y} \in \Sigma \setminus \{\bar{x}\}$ satisfying $Z_\kappa(\bar{x},\bar{y}) = 0$. Therefore, Proposition \ref{key.estimate} implies that 
\begin{align*} 
0 &\leq \sum_{i=1}^2 \frac{\partial^2 Z_\kappa}{\partial x_i^2}(\bar{x},\bar{y}) + 2 \sum_{i=1}^2 \frac{\partial^2 Z_\kappa}{\partial x_i \, \partial y_i}(\bar{x},\bar{y}) + \sum_{i=1}^2 \frac{\partial^2 Z_\kappa}{\partial y_i^2}(\bar{x},\bar{y}) \\ 
&\leq -\frac{\kappa^2-1}{\kappa} \, \frac{\Psi(\bar{x})}{1-\langle F(\bar{x}),F(\bar{y}) \rangle} \sum_{i=1}^2 \Big \langle \frac{\partial F}{\partial x_i}(\bar{x}),F(\bar{y}) \Big \rangle^2, 
\end{align*} 
where $(x_1,x_2)$ and $(y_1,y_2)$ are suitable coordinate systems around $\bar{x}$ and $\bar{y}$, respectively. This gives 
\[\Big \langle \frac{\partial F}{\partial x_i}(\bar{x}),F(\bar{y}) \Big \rangle = 0\] 
for $i = 1,2$. Using (\ref{gradient.a}), it follows that $\frac{\partial \Psi}{\partial x_i}(\bar{x}) = 0$ for $i = 1,2$. Thus, the gradient of $\Psi$ vanishes at each point in $\Omega$. By the unique continuation theorem for elliptic partial differential equations (cf. \cite{Aronszajn}), the gradient of $\Psi$ vanishes identically. From this, we deduce that the surface is congruent to the Clifford torus. This completes the proof of Theorem \ref{lawson.conj}. \\

The proof of the Lawson conjecture can be extended to give a classification of all Alexandrov immersed minimal tori in $S^3$:

\begin{theorem}[S.~Brendle \cite{Brendle-note}]
\label{extension.of.Lawson.conj}
Let $F: \Sigma \to S^3$ be an immersed minimal surface in $S^3$ of genus $1$. Moreover, we assume that $F$ is an Alexandrov immersion in the sense of Definition \ref{definition.of.alexandrov.immersion} above. Then $\Sigma$ is rotationally symmetric.
\end{theorem} 

In the remainder of this section, we will describe the proof of Theorem \ref{extension.of.Lawson.conj}. As usual, we will identify $S^3$ with the unit sphere in $\mathbb{R}^4$. By assumption, there exists a compact manifold $N$ and an immersion $\bar{F}: N \to S^3$ such that $\partial N = \Sigma$ and $\bar{F}|_\Sigma = F$. It will be convenient to put a Riemannian metric on $N$ so that $\bar{F}$ is a local isometry. Since $F$ is a local isometry, we can find a real number $\delta>0$ so that $\bar{F}(x) \neq \bar{F}(y)$ for all points $x,y \in N$ satisfying $d_N(x,y) \in (0,\delta)$.

For each point $x \in \Sigma$, we denote by $\nu(x) \in T_{F(x)} S^3$ the push-forward of the outward-pointing unit normal to $\Sigma$ at the point $x$ under the map $\bar{F}$. Given any point $x \in \Sigma$ and any number $\alpha \geq 1$, we define 
\[D_\alpha(x) = \big \{ p \in S^3: \alpha \, \Psi(x) \, (1 - \langle F(x),p \rangle) + \langle \nu(x),p \rangle \leq 0 \big \}.\] 
Note that $D_\alpha(x)$ is a closed geodesic ball in $S^3$ with radius less than $\frac{\pi}{2}$. Moreover, the point $F(x)$ lies on the boundary $\partial D_\alpha(x)$, and the outward-pointing unit normal vector to $\partial D_\alpha(x)$ at the point $F(x)$ is given by $\nu(x)$. 

Let $I$ denote the set of all points $(x,\alpha) \in \Sigma \times [1,\infty)$ with the property that there exists a smooth map $G: D_\alpha(x) \to N$ such that $\bar{F} \circ G = \text{\rm id}_{D_\alpha(x)}$ and $G(F(x)) = x$. 

\begin{lemma}
\label{aux.1}
Let us fix a pair $(x,\alpha) \in I$. Then there is a unique map $G: D_\alpha(x) \to N$ such that $\bar{F} \circ G = \text{\rm id}_{D_\alpha(x)}$ and $G(F(x)) = x$.
\end{lemma}

\textbf{Proof of Lemma \ref{aux.1}.} 
It suffices to prove the uniqueness statement. Suppose that $G$ and $\tilde{G}$ are two maps which have the required properties. Then $\bar{F}(G(p)) = \bar{F}(\tilde{G}(p)) = p$ for all points $p \in D_\alpha(x)$. This implies $d_N(G(p),\tilde{G}(p)) \notin (0,\delta)$ for all $p \in D_\alpha(x)$. By continuity, we either have $G(p) = \tilde{G}(p)$ for all $p \in D_\alpha(x)$ or we have $G(p) \neq \tilde{G}(p)$ for all $p \in D_\alpha(x)$. Since $G(F(x)) = \tilde{G}(F(x)) = x$, the second case cannot occur. Thus, we conclude that $G(p) = \tilde{G}(p)$ for all $p \in D_\alpha(x)$. \\

\begin{lemma}
\label{aux.2}
The set $I$ is closed. Moreover, the map $G$ depends continuously on the pair $(x,\alpha)$.
\end{lemma}

\textbf{Proof of Lemma \ref{aux.2}.} 
Let us consider a sequence of pairs $(x^{(m)},\alpha^{(m)}) \in I$ such that $\lim_{m \to \infty} (x^{(m)},\alpha^{(m)}) = (\bar{x},\bar{\alpha})$. For each $m$, we can find a smooth map $G^{(m)}: D_{\alpha^{(m)}}(x^{(m)}) \to N$ such that $\bar{F} \circ G^{(m)} = \text{\rm id}_{D_{\alpha^{(m)}}(x^{(m)})}$ and $G^{(m)}(F(x^{(m)})) = x^{(m)}$. Since $\bar{F}$ is a smooth immersion, the maps $G^{(m)}$ are uniformly bounded in $C^2$ norm. Hence, after passing to a subsequence, the maps $G^{(m)}$ converge in $C^1$ to a map $G: D_{\bar{\alpha}}(\bar{x}) \to N$ satisfying $\bar{F} \circ G = \text{\rm id}_{D_{\bar{\alpha}}(\bar{x})}$ and $G(F(\bar{x})) = \bar{x}$. It is easy to see that the map $G$ is smooth. Thus, $(\bar{x},\bar{\alpha}) \in I$, and the assertion follows. \\

In the next step, we show that the set $I$ is non-empty. \\

\begin{lemma}
\label{aux.3}
We have $(x,\alpha) \in I$ if $\alpha$ is sufficiently large. 
\end{lemma}

\textbf{Proof of Lemma \ref{aux.3}.} 
By Proposition \ref{absence.of.umbilic.points}, the function $\Psi$ is strictly positive. Hence, the radius of the geodesic ball $D_\alpha(x) \subset S^3$ will be arbitrarily small if $\alpha$ is sufficiently large. Hence, if $\alpha$ is large enough, we can use the implicit function theorem to construct a smooth map $G: D_\alpha(x) \to N$ such that $\bar{F} \circ G = \text{\rm id}_{D_\alpha(x)}$ and $G(F(x)) = x$. This proves Lemma \ref{aux.3}. \\

We now continue with the proof of Theorem \ref{extension.of.Lawson.conj}. Let  
\[\kappa = \inf \{\alpha: \text{\rm $(x,\alpha) \in I$ for all $x \in \Sigma$}\}.\] 
Clearly, $\kappa \in [1,\infty)$. For each point $x \in \Sigma$, there is a unique map $G_x: D_\kappa(x) \to N$ such that $\bar{F} \circ G_x = \text{\rm id}_{D_\kappa(x)}$ and $G_x(F(x)) = x$. For each point $x \in \Sigma$, the map $G_x$ and the map $\bar{F}|_{G_x(D_\kappa(x))}$ are injective. To complete the proof, we distinguish two cases: \\

\textit{Case 1:} We first consider the special case that $\kappa=1$. We begin with a lemma: 

\begin{lemma} 
\label{Z.nonnegative.near.diagonal}
Given any point $\bar{x} \in \Sigma$, there exists an open set $V$ containing $\bar{x}$ such that $Z_1(\bar{x},y) \geq 0$ for all $y \in V$.
\end{lemma}

\textbf{Proof of Lemma \ref{Z.nonnegative.near.diagonal}.} 
We argue by contradiction. Suppose that there exists a sequence of points $y^{(m)} \in \Sigma$ such that $\lim_{m \to \infty} y^{(m)} = \bar{x}$ and $Z_1(\bar{x},y^{(m)}) < 0$ for all $m$. Since $Z_1(\bar{x},y^{(m)}) < 0$, the point $F(y^{(m)})$ lies in the interior of the geodesic ball $D_1(\bar{x})$. Therefore, the point $\tilde{y}^{(m)} := G_{\bar{x}}(F(y^{(m)}))$ lies in the interior of $N$. Since $y^{(m)}$ lies on the boundary $\partial N = \Sigma$, it follows that 
\[\tilde{y}^{(m)} \neq y^{(m)}.\] 
On the other hand, we have 
\[\bar{F}(\tilde{y}^{(m)}) = F(y^{(m)})\] 
and 
\[\lim_{m \to \infty} \tilde{y}^{(m)} = \lim_{m \to \infty} G_{\bar{x}}(F(y^{(m)})) = G_{\bar{x}}(F(\bar{x})) = \bar{x} = \lim_{m \to \infty} y^{(m)}.\] 
This contradicts the fact that $\bar{F}$ is an immersion. \\

\begin{lemma} 
\label{rotational.symmetry}
Fix a point $\bar{x} \in \Sigma$, and let $\{e_1,e_2\}$ is an orthonormal basis of $T_{\bar{x}} \Sigma$ such that $h(e_1,e_1) = \Psi(\bar{x})$, $h(e_1,e_2) = 0$, and $h(e_2,e_2) = -\Psi(\bar{x})$. Then $\langle \nabla \Psi(\bar{x}),e_1 \rangle = 0$.
\end{lemma}

\textbf{Proof of Lemma \ref{rotational.symmetry}.} 
For abbreviation, we define a vector $\xi \in \mathbb{R}^4$ by 
\[\xi = \Psi(\bar{x}) \, F(\bar{x}) - \nu(\bar{x}).\] 
Note that $\xi$ is orthogonal to the tangent plane to $dF_{\bar{x}}(e_1)$ and $dF_{\bar{x}}(e_2)$. Let $\sigma: \mathbb{R} \to \Sigma$ be a geodesic such that $\sigma(0) = \bar{x}$ and $\sigma'(0) = e_1$. By Lemma \ref{Z.nonnegative.near.diagonal}, we have $Z_1(\bar{x},y) \geq 0$ if $y$ is sufficiently close to $\bar{x}$. Consequently, the function 
\[f(t) = Z_1(\sigma(0),\sigma(t)) = \Psi(\bar{x}) - \langle \xi,F(\sigma(t)) \rangle\] 
is nonnegative when $t$ is sufficiently small. As above, we compute 
\[f'(t) = -\langle \xi,dF_{\sigma(t)}(\sigma'(t)) \rangle,\]
\[f''(t) = \langle \xi,F(\sigma(t)) \rangle + h(\sigma'(t),\sigma'(t)) \, \langle \xi,\nu(\sigma(t)) \rangle,\] 
and 
\begin{align*} 
f'''(t) 
&= \langle \xi,dF_{\sigma(t)}(\sigma'(t)) \rangle + h(\sigma'(t),\sigma'(t)) \, \langle \xi,D_{\sigma'(t)} \nu \rangle \\ 
&+ (D_{\sigma'(t)}^\Sigma h)(\sigma'(t),\sigma'(t)) \, \langle \xi,\nu(\sigma(t)) \rangle. 
\end{align*} 
Setting $t=0$, we conclude that $f(0) = f'(0) = f''(0) = 0$. Since the function $f(t)$ is nonnegative in a neighborhood of $0$, it follows that $f'''(0) = 0$. This implies that $(D_{e_1}^\Sigma h)(e_1,e_1) = 0$. From this, the assertion follows easily. \\

Using Lemma \ref{rotational.symmetry}, we conclude that the function $\Psi$ is constant along one set of curvature lines on $\Sigma$. This implies that $\Sigma$ is rotationally symmetric. \\

\textit{Case 2:} We next consider the case $\kappa > 1$. In order to handle this case, we need several auxiliary results:

\begin{lemma} 
\label{technical.ingredient}
There exists a constant $\beta > 0$ with the following property: if $x \in \Sigma$ and $p \in \partial D_\kappa(x)$ are two points satisfying $|p - F(x)| \leq \beta$, then we have $d_N(G_x(p),\Sigma) \geq \beta \, |p - F(x)|^2$.
\end{lemma}

\textbf{Proof of Lemma \ref{technical.ingredient}.} 
Let us fix a point $\bar{x} \in \Sigma$. We consider the function 
\[\rho: \partial D_\kappa(\bar{x}) \to \mathbb{R}, \quad p \mapsto d_N(G_{\bar{x}}(p),\Sigma).\] 
Clearly, $\rho(F(\bar{x})) = 0$, and the gradient of the function $\rho$ at the point $F(\bar{x})$ vanishes. Moreover, since $\kappa > 1$, the Hessian of the function $\rho$ at the point $F(\bar{x})$ is positive definite. Hence, we can find a positive constant $\beta > 0$ such that $\rho(p) \geq \beta \, |p - F(\bar{x})|^2$ for all points $p \in \partial D_\kappa(\bar{x})$ satisfying $|p - F(\bar{x})| \leq \beta$. This completes the proof of Lemma \ref{technical.ingredient}. \\

\begin{lemma} 
\label{touch}
There exists a point $\hat{x} \in \Sigma$ such that $\Sigma \cap G_{\hat{x}}(\partial D_\kappa(\hat{x})) \neq \{\hat{x}\}$.
\end{lemma}

\textbf{Proof of Lemma \ref{touch}.} 
Suppose this is false. Then $\Sigma \cap G_x(\partial D_\kappa(x)) = \{x\}$ for all $x \in \Sigma$. This implies that $d_N(G_x(p),\Sigma) > 0$ for all $x \in \Sigma$ and all points $p \in \partial D_\kappa(x) \setminus \{F(x)\}$. Using Lemma \ref{technical.ingredient}, we conclude that there exists a positive constant $\gamma > 0$ such that $d_N(G_x(p),\Sigma) \geq \gamma \, |p - F(x)|^2$ for all points $x \in \Sigma$ and all points $p \in \partial D_\kappa(x)$. Hence, if $\varepsilon > 0$ is sufficiently small, then the map $G_x: D_\kappa(x) \to N$ can be extended to a smooth map $\tilde{G}_x: D_{\kappa-\varepsilon}(x) \to N$ satisfying $\bar{F} \circ \tilde{G}_x = \text{\rm id}_{D_{\kappa-\varepsilon}(x)}$. Consequently, $(x,\kappa-\varepsilon) \in I$ for all $x \in \Sigma$. This contradicts the definition of $\kappa$. \\

Let $\hat{x} \in \Sigma$ be chosen as in Lemma \ref{touch}. Moreover, let us pick a point $\hat{y} \in \Sigma \cap G_{\hat{x}}(\partial D_\kappa(\hat{x}))$ such that $\hat{x} \neq \hat{y}$. Since $\hat{y} \in G_{\hat{x}}(\partial D_\kappa(\hat{x}))$, we conclude that $F(\hat{y}) \in \partial D_\kappa(\hat{x})$ and $G_{\hat{x}}(F(\hat{y})) = \hat{y}$. Moreover, we claim that $F(\hat{x}) \neq F(\hat{y})$; indeed, if $F(\hat{x}) = F(\hat{y})$, then $\hat{x} = G_{\hat{x}}(F(\hat{x})) = G_{\hat{x}}(F(\hat{y})) = \hat{y}$, which contradicts our choice of $\hat{y}$.

\begin{lemma}
\label{Z.nonnegative}
We can find open sets $U,V \subset \Sigma$ such that $\hat{x} \in U$, $\hat{y} \in V$, and $Z_\kappa(x,y) \geq 0$ for all points $(x,y) \in U \times V$.
\end{lemma}

\textbf{Proof of Lemma \ref{Z.nonnegative}.} 
We argue by contradiction. Suppose that there exist sequences of points $x^{(m)},y^{(m)} \in \Sigma$ such that $\lim_{m \to \infty} x^{(m)} = \hat{x}$, $\lim_{m \to \infty} y^{(m)} = \hat{y}$, and $Z_\kappa(x^{(m)},y^{(m)}) < 0$. Since $Z_\kappa(x^{(m)},y^{(m)}) < 0$, the point $F(y^{(m)})$ lies in the interior of the geodesic ball $D_\kappa(x^{(m)})$. Therefore, the point $\tilde{y}^{(m)} := G_{x^{(m)}}(F(y^{(m)}))$ lies in the interior of $N$. Since the point $y^{(m)}$ lies on the boundary $\partial N = \Sigma$, we conclude that 
\[\tilde{y}^{(m)} \neq y^{(m)}.\] 
On the other hand, we have 
\[\bar{F}(\tilde{y}^{(m)}) = F(y^{(m)})\] 
and 
\[\lim_{m \to \infty} \tilde{y}^{(m)} = \lim_{m \to \infty} G_{x^{(m)}}(F(y^{(m)})) = G_{\hat{x}}(F(\hat{y})) = \hat{y} = \lim_{m \to \infty} y^{(m)}\] 
by Lemma \ref{aux.2}. This contradicts the fact that $\bar{F}$ is an immersion. Thus, $Z_\kappa(x,y) \geq 0$ if $(x,y)$ is sufficiently close to $(\hat{x},\hat{y})$. This completes the proof of Lemma \ref{Z.nonnegative}. \\

Since $F(\hat{x}) \neq F(\hat{y})$, we can choose the sets $U$ and $V$ small enough so that $F(\bar{U}) \cap F(\bar{V}) = \emptyset$. We now define 
\[\Omega = \{x \in U: \text{\rm there exists a point $y \in V$ such that $Z_\kappa(x,y) = 0$}\}.\] 
Since $F(\hat{y}) \in \partial D_\kappa(\hat{x})$, we have $Z_\kappa(\hat{x},\hat{y}) = 0$. Consequently,  $\hat{x} \in \Omega$. In particular, the set $\Omega$ is non-empty. Using Proposition \ref{key.estimate} and Bony's version of the strict maximum principle (cf. \cite{Bony}), we conclude that the set $\Omega$ is open. 

As above, we will show that the gradient of $\Psi$ vanishes at each point $\bar{x} \in \Omega$. To see this, we consider a pair of points $\bar{x} \in U$ and $\bar{y} \in V$ satisfying $Z_\kappa(\bar{x},\bar{y}) = 0$. Using Proposition \ref{key.estimate}, we obtain 
\begin{align*} 
0 &\leq \sum_{i=1}^2 \frac{\partial^2 Z_\kappa}{\partial x_i^2}(\bar{x},\bar{y}) + 2 \sum_{i=1}^2 \frac{\partial^2 Z_\kappa}{\partial x_i \, \partial y_i}(\bar{x},\bar{y}) + \sum_{i=1}^2 \frac{\partial^2 Z_\kappa}{\partial y_i^2}(\bar{x},\bar{y}) \\ 
&\leq -\frac{\kappa^2-1}{\kappa} \, \frac{\Psi(\bar{x})}{1-\langle F(\bar{x}),F(\bar{y}) \rangle} \sum_{i=1}^2 \Big \langle \frac{\partial F}{\partial x_i}(\bar{x}),F(\bar{y}) \Big \rangle^2, 
\end{align*} 
where $(x_1,x_2)$ and $(y_1,y_2)$ are suitable coordinate systems around $\bar{x}$ and $\bar{y}$, respectively. From this, we deduce that 
\[\Big \langle \frac{\partial F}{\partial x_i}(\bar{x}),F(\bar{y}) \Big \rangle = 0\] 
for $i = 1,2$. Using (\ref{gradient.a}), we conclude that $\nabla \Psi(\bar{x}) = 0$ for each point $\bar{x} \in \Omega$. Hence, it follows from standard unique continuation arguments (cf. \cite{Aronszajn}) that the gradient of $\Psi$ vanishes identically. This implies that $F$ is congruent to the Clifford torus. This completes the proof of Theorem \ref{extension.of.Lawson.conj}. \\

We note that all the results in this section have analogues for surfaces with constant mean curvature. For example, the proof of Almgren's theorem (Theorem \ref{almgren.genus.0}) can be adapted to show that an immersed constant mean curvature surface in $S^3$ of genus $0$ is a geodesic sphere. Similarly, there is a generalization of Proposition \ref{absence.of.umbilic.points} which asserts that a constant mean curvature surface of genus $1$ has no umbilic points, and the norm of the trace-free part of the second fundamental form still satisfies a Simons-type identity. Andrews and Li \cite{Andrews-Li} observed that the proof of Theorem \ref{lawson.conj} can be adapted to show that any embedded constant mean curvature surface of genus $1$ is rotationally symmetric. More generally, it was shown in \cite{Brendle-note} that any Alexandrov immersed constant mean curvature surface in $S^3$ is rotationally symmetric. 

Finally, in a recent paper \cite{Brendle1}, we obtained a uniqueness theorem for embedded constant mean curvature surfaces in certain rotationally symmetric spaces. This result generalizes the classical Alexandrov theorem in Euclidean space. Moreover, there is a rich literature on constant mean curvature surfaces in asymptotically flat three-manifolds; see e.g. \cite{Brendle-Eichmair}, \cite{Christodoulou-Yau}, \cite{Eichmair-Metzger1}, \cite{Eichmair-Metzger2}, \cite{Huisken-Yau}, \cite{Qing-Tian}.

\section{Estimates for the Morse index and area of a minimal surface and the Willmore conjecture} 

The Willmore energy of a two-dimensional surface $\Sigma$ in $S^3$ is defined by 
\begin{equation} 
\label{willmore.1}
\mathscr{W}(\Sigma) = \int_\Sigma \Big ( 1 + \frac{H^2}{4} \Big ), 
\end{equation}
where $H$ denotes the mean curvature of $\Sigma$. Note that $\mathscr{W}(\Sigma) = 4\pi$ for the equator, and $\mathscr{W}(\Sigma) = 2\pi^2$ for the Clifford torus. 

We first collect some classical facts about the Willmore functional. The Gauss equations imply that 
\[1 + \frac{H^2}{4} = K + \frac{|\mathring{A}|^2}{2},\] 
where $K$ is the intrinsic Gaussian curvature of $\Sigma$ and $\mathring{A}$ denotes the trace-free part of the second fundamental form. Thus, 
\begin{equation} 
\label{willmore.2}
\mathscr{W}(\Sigma) = 2\pi \chi(\Sigma) + \int_\Sigma \frac{|\mathring{A}|^2}{2} 
\end{equation}
by the Gauss-Bonnet theorem. The identity (\ref{willmore.2}) shows that the Willmore functional is invariant under conformal transformations in $S^3$. More precisely, let us consider a conformal transformation $\psi: S^3 \to S^3$ of the form 
\[\psi(x) = a + \frac{1-|a|^2}{1+2 \, \langle a,x \rangle+|a|^2} \, (x+a),\] 
where $a$ is a vector $a \in \mathbb{R}^4$ satisfying $|a| < 1$. Then 
\[\mathscr{W}(\psi(\Sigma)) = \mathscr{W}(\Sigma)\] 
for any surface $\Sigma \subset S^3$.

The following result is well-known (see e.g. \cite{Li-Yau}):

\begin{proposition}
\label{lower.bound}
Let $\Sigma$ be an immersed surface in $S^3$, and let $p$ be a point on $\Sigma$. Then $\mathscr{W}(\Sigma) \geq 4\pi m$, where $m$ is the multiplicity of $\Sigma$ at $p$. In particular, $\mathscr{W}(\Sigma) \geq 4\pi$. Moreover, if $\mathscr{W}(\Sigma) < 8\pi$, then $\Sigma$ is embedded.
\end{proposition}

\textbf{Sketch of the proof of Proposition \ref{lower.bound}.} 
Let $\psi: S^3 \to S^3$ be a conformal transformation of the form 
\[\psi(x) = a + \frac{1-|a|^2}{1+2 \, \langle a,x \rangle+|a|^2} \, (x+a),\] 
where $a$ is a vector $a \in \mathbb{R}^4$ satisfying $|a| < 1$. The conformal invariance of the Willmore functional implies 
\[\int_\Sigma \Big ( \frac{1-|a|^2}{1+2 \, \langle a,x \rangle+|a|^2} \Big )^2 = \text{\rm area}(\psi(\Sigma)) \leq \mathscr{W}(\psi(\Sigma)) = \mathscr{W}(\Sigma).\] 
If we put $a = -(1-\varepsilon) \, p$ and take the limit as $\varepsilon \to 0$, we conclude that $4\pi m \leq \mathscr{W}(\Sigma)$, as claimed. \\

In \cite{Ros2}, Ros discovered a connection between the Willmore energy of a surface $\Sigma$ and the area of a distance surface: 

\begin{proposition}[A.~Ros \cite{Ros2}]
\label{distance.surface}
Let $\Sigma$ be an immersed surface in $S^3$. Moreover, let $\nu(x)$ be the unit normal vector field along $\Sigma$, and let
\[\Sigma_t = \{\cos t \, x + \sin t \, \nu(x): x \in \Sigma\}.\] 
Then 
\[\text{\rm area}(\Sigma_t) \leq \mathscr{W}(\Sigma)\] 
for $t \in (-\pi,\pi)$.
\end{proposition} 

\textbf{Sketch of the proof of Proposition \ref{distance.surface}.}
The area of $\Sigma_t$ is given by 
\[\text{\rm area}(\Sigma_t) = \int_\Sigma (\cos t + \sin t \, \lambda_1) \, (\cos t + \sin t \, \lambda_2),\] 
where $\lambda_1$ and $\lambda_2$ denote the principal curvatures of $\Sigma$. We next compute 
\begin{align*} 
&(\cos t + \sin t \, \lambda_1) \, (\cos t + \sin t \, \lambda_2) \\ 
&= 1 + \Big ( \frac{\lambda_1+\lambda_2}{2} \Big )^2 - \sin^2 t \, \Big ( \frac{\lambda_1-\lambda_2}{2} \Big )^2 \\ 
&- \Big ( \sin t - \cos t \, \frac{\lambda_1+\lambda_2}{2} \Big )^2 \\ 
&\leq 1 + \frac{H^2}{4}. 
\end{align*} 
Thus, 
\[\text{\rm area}(\Sigma_t) \leq \int_\Sigma \Big ( 1 + \frac{H^2}{4} \Big ) = \mathscr{W}(\Sigma),\] 
as claimed. \\

Combining Theorem \ref{distance.surface} with the solution of the isoperimetric problem in $\mathbb{RP}^3$ in \cite{Ritore-Ros}, Ros was able to give a sharp lower bound for the Willmore energy when $\Sigma$ has antipodal symmetry: 

\begin{theorem}[A.~Ros \cite{Ros2}]
\label{willmore.antipodal.symmetry}
Suppose that $\Sigma$ is an embedded surface of genus $1$ which is invariant under antipodal reflection. Then $\mathscr{W}(\Sigma) \geq 2\pi^2$. 
\end{theorem}

\textbf{Sketch of the proof of Theorem \ref{willmore.antipodal.symmetry}.}
The surface $\Sigma$ divides $S^3$ into two regions, which we denote by $N$ and $\tilde{N}$. Since the genus of $\Sigma$ is odd, the quotient of $\Sigma$ under the natural $\mathbb{Z}_2$ action is an orientable surface in $\mathbb{RP}^3$. Hence, there exists a unit normal vector field $\nu$ along $\Sigma$ which is invariant under antipodal reflection. Consequently, both $N$ and $\tilde{N}$ are invariant under antipodal reflection.

Without loss of generality, we may assume that $\text{\rm vol}(N) \leq \frac{1}{2} \, \text{\rm vol}(S^3)$. Hence, we can find a real number $t \in [0,\pi)$ such that $\text{\rm vol}(N_t) = \frac{1}{2} \, \text{\rm vol}(S^3)$, where 
\[N_t = \{x \in S^3: d(x,N) \leq t\}.\] 
By a theorem of Ritor\'e and Ros, any region in $S^3$ which has volume $\frac{1}{2} \, \text{\rm vol}(S^3)$ and is invariant under antipodal symmetry has boundary area at least $2\pi^2$ (see \cite{Ritore-Ros} or \cite{Ros3}, Corollary 5). Therefore, $\text{\rm area}(\partial N_t) \geq 2\pi^2$. On the other hand, the boundary of $N_t$ is contained in the set 
\[\Sigma_t = \{\cos t \, x + \sin t \, \nu(x): x \in \Sigma\}.\] 
Using Proposition \ref{distance.surface}, we conclude that 
\[\mathscr{W}(\Sigma) \geq \text{\rm area}(\Sigma_t) \geq \text{\rm area}(\partial N_t) \geq 2\pi^2,\] 
as claimed. \\

In 1965, Willmore proposed the problem of minimizing the Willmore energy among surfaces of genus $1$. This led him to the following conjecture:

\begin{conjecture}[T.J.~Willmore \cite{Willmore1}, \cite{Willmore2}]
Let $\Sigma$ be a surface in $S^3$ with genus $1$. Then $\mathscr{W}(\Sigma) \geq 2\pi^2$.
\end{conjecture} 

Theorem \ref{willmore.antipodal.symmetry} shows that the Willmore conjecture holds for tori with antipodal symmetry. We note that Topping \cite{Topping} has obtained an alternative proof of Theorem \ref{willmore.antipodal.symmetry}, which is based on techniques from integral geometry. In 2012, Marques and Neves \cite{Marques-Neves} verified the Willmore conjecture in full generality. Their proof relies on the min-max theory for minimal surfaces. The argument in \cite{Marques-Neves} also uses a sharp estimate for the Morse index of a minimal surface in $S^3$, which we describe below. Recall that the Jacobi operator of a minimal surface in $S^3$ is defined by $L = -\Delta_\Sigma - |A|^2 - 2$. Moreover, the Morse index of a minimal surface is defined as the number of negative eigenvalues of the Jacobi operator, counted according to multiplicity.

The following theorem, due to Urbano, characterizes the Clifford torus as the unique minimal surface in $S^3$ which has genus at least $1$ and Morse index at most $5$. 

\begin{theorem}[F.~Urbano \cite{Urbano}]
\label{index.estimate}
Let $\Sigma$ be an immersed minimal surface in $S^3$ of genus at least $1$. Then the Morse index of $\Sigma$ is at least $5$. Moreover, the Morse index of $\Sigma$ is equal to $5$ if and only if $\Sigma$ is congruent to the Clifford torus.
\end{theorem}

\textbf{Proof of Theorem \ref{index.estimate}.} 
Let $U \subset C^\infty(\Sigma)$ be the space of all functions of the form $\langle a,\nu \rangle$, where $a$ is a fixed vector in $\mathbb{R}^4$ and $\nu$ denotes the unit normal vector to $\Sigma$. Since $\Sigma$ is not totally geodesic, we have $\dim U = 4$. Moreover, every function $u \in U$ satisfies $\Delta_\Sigma u + |A|^2 \, u = 0$, hence $Lu = -2u$. Thus, $-2$ is an eigenvalue of the Jacobi operator $L$, and the associated eigenspace has dimension at least $4$. However, the first eigenvalue $\lambda_1$ of $L$ has multiplicity $1$. Therefore, $\lambda_1 < -2$, and $L$ has at least five negative eigenvalues. 

Suppose now that the Jacobi operator $L$ has exactly five negative eigenvalues. Let $\rho$ denote the eigenfunction associated with the eigenvalue $\lambda_1$. Note that $\rho$ is a positive function. We consider a conformal transformation $\psi: S^3 \to S^3$ of the form 
\[\psi(x) = a + \frac{1-|a|^2}{1+2 \, \langle a,x \rangle+|a|^2} \, (x+a),\] 
where $a$ is a vector $a \in \mathbb{R}^4$ satisfying $|a| < 1$. We can choose the vector $a$ in a such a way that 
\[\int_\Sigma \rho \, \psi_i(x) = 0\] 
for $i \in \{1,2,3,4\}$, where $\psi_i(x)$ denotes the $i$-th component of the vector $\psi(x) \in S^3 \subset \mathbb{R}^4$. 

By assumption, $L$ has exactly five negative eigenvalues. In particular, $L$ has no eigenvalues between $\lambda_1$ and $2$. Since the function $\psi_i$ is orthogonal to the eigenfunction $\rho$, we conclude that  
\begin{equation} 
\label{stab}
\int_\Sigma (|\nabla^\Sigma \psi_i|^2 - |A|^2 \, \psi_i^2) = \int_\Sigma \psi_i \, (L\psi_i + 2\psi_i) \geq 0 
\end{equation}
for each $i \in \{1,2,3,4\}$. On the other hand, the conformal invariance of the Willmore functional implies that 
\begin{equation} 
\label{sum.1}
\sum_{i=1}^4 \int_\Sigma |\nabla^\Sigma \psi_i|^2 = 2 \, \text{\rm area}(\psi(\Sigma)) \leq 2 \, \mathscr{W}(\psi(\Sigma)) = 2 \, \mathscr{W}(\Sigma) = 2 \, \text{\rm area}(\Sigma). 
\end{equation}
Moreover, it follows from the Gauss-Bonnet theorem that 
\begin{equation} 
\label{sum.2}
\sum_{i=1}^4 \int_\Sigma |A|^2 \, \psi_i^2 = \int_\Sigma |A|^2 = 2 \int_\Sigma (1-K) \geq 2 \, \text{\rm area}(\Sigma). 
\end{equation}
Combining the inequalities (\ref{sum.1}) and (\ref{sum.2}) gives 
\begin{equation} 
\label{sum.3}
\sum_{i=1}^4 \int_\Sigma (|\nabla^\Sigma \psi_i|^2 - |A|^2 \, \psi_i^2) \leq 0. 
\end{equation}
Putting these facts together, we conclude that all the inequalities must, in fact, be equalities. In particular, we must have $\mathscr{W}(\psi(\Sigma)) = \text{\rm area}(\psi(\Sigma))$. Consequently, the surface $\psi(\Sigma)$ must have zero mean curvature. This implies that $\langle a,\nu \rangle = 0$ at each point on $\Sigma$. Since $\Sigma$ is not totally geodesic, it follows that $a=0$. Furthermore, since $\int_\Sigma \rho \, \psi_i = 0$ and $\int_\Sigma (|\nabla^\Sigma \psi_i|^2 - |A|^2 \, \psi_i^2) = 0$, we conclude that the function $\psi_i$ is an eigenfunction of the Jacobi operator with eigenvalue $-2$. Consequently, $\Delta_\Sigma \psi_i + |A|^2 \, \psi_i = 0$ for each $i \in \{1,2,3,4\}$. Since $a=0$, we conclude that $\Delta_\Sigma x_i + |A|^2 \, x_i = 0$. Since $\Delta_\Sigma x_i + 2 \, x_i = 0$, we conclude that $|A|^2=2$ and the Gaussian curvature of $\Sigma$ vanishes. This implies that $\Sigma$ is the Clifford torus. 

Finally, it is straightforward to verify that the Jacobi operator on the Clifford torus has exactly five negative eigenvalues. This completes the proof of Theorem \ref{index.estimate}. \\

Theorem \ref{index.estimate} gives a lower bound for the number of negative eigenvalues of the Jacobi operator $L = -\Delta_\Sigma - |A|^2 - 2$. It is an interesting problem to understand the nullspace of $L$. Clearly, if $K$ is an ambient rotation vector field, then the function $\langle K,\nu \rangle$ lies in the nullspace of $L$. It is a natural to conjecture that the nullspace of $L$ should consist precisely of the functions $\langle K,\nu \rangle$ where $K$ is an ambient rotation vector field.

We now describe the min-max procedure of Marques and Neves \cite{Marques-Neves}. Let us fix an embedded surface $\Sigma$ in $S^3$, and let $\nu$ be the unit normal vector field along $\Sigma$. Given any point $a$ in the open unit ball $B^4$, we consider the conformal transformation 
\[\psi(x) = a + \frac{1-|a|^2}{1+2 \, \langle a,x \rangle+|a|^2} \, (x+a).\] 
For each $t \in (-\pi,\pi)$, we denote by $\Sigma_{(a,t)}$ the parallel surface to $\psi(\Sigma)$ at distance $t$. This defines a five-parameter family of surfaces in $S^3$, which is parametrized by $B^4 \times (-\pi,\pi)$. However, the map $(a,t) \mapsto \Sigma_{(a,t)}$ does not extend continuously to $\bar{B}^4 \times [-\pi,\pi]$.

In order to overcome this obstacle, Marques and Neves consider the map 
\[\Psi: \Sigma \times [0,1] \times [-\pi,\pi] \to \bar{B}^4, \quad (x,r,s) \mapsto (1-r) \, (\cos s \, x + \sin s \, \nu(x)).\] 
For abbreviation, let 
\[\Omega_\varepsilon = \{\Psi(x,r,s): x \in \Sigma, \, r \geq 0, \, \sqrt{r^2+s^2} \leq \varepsilon\}.\] 
Moreover, let $T: \bar{B}^4 \to \bar{B}^4$ be a continuous map with the following properties: 
\begin{itemize}
\item $T = \text{\rm id}$ on $\bar{B}^4 \setminus \Omega_{2\varepsilon}$.
\item $T$ maps the point $\Psi(x,r,s) \in \Omega_{2\varepsilon} \setminus \Omega_\varepsilon$ to the point $\Psi(x,\tilde{r},\tilde{s}) \in \Omega_{2\varepsilon}$, where $\tilde{r} = (\frac{\sqrt{r^2+s^2}}{\varepsilon} - 1) \, r$ and $\tilde{s} = (\frac{\sqrt{r^2+s^2}}{\varepsilon} - 1) \, s$.
\item $T$ maps the point $\Psi(x,r,s) \in \Omega_\varepsilon$ to the point $x \in \Sigma$.
\end{itemize}
Given any pair $(a,t) \in B^4 \times (-\pi,\pi)$, one can define a surface $\hat{\Sigma}_{(a,t)}$ in the following way: 
\begin{itemize} 
\item Suppose first that $a \in B^4 \setminus \Omega_\varepsilon$. In this case, one defines $\hat{\Sigma}_{(a,t)} = \Sigma_{(T(a),t)}$. 
\item Suppose next that $a \in \Omega_\varepsilon$. Let $a = \Psi(x,r,s)$, where $x \in \Sigma$ and $\sqrt{r^2+s^2} \leq \varepsilon$. In this case, one defines $\hat{\Sigma}_{(a,t)}$ to be a geodesic sphere centered at the point $-\sin \theta \, x - \cos \theta \, \nu(x)$ of radius $\frac{\pi}{2} - \theta + t$, where $\theta = \arcsin(\frac{s}{\varepsilon}) \in [-\frac{\pi}{2},\frac{\pi}{2}]$.
\end{itemize}
The family of surfaces $\hat{\Sigma}_{(a,t)}$ is called the canonical family associated with $\Sigma$. Its main properties are summarized in the following proposition: 

\begin{proposition}[F.C.~Marques, A.~Neves \cite{Marques-Neves}]
\label{canonical.family}
The canonical family has the following properties: 
\begin{itemize}
\item The map $(a,t) \mapsto \hat{\Sigma}_{(a,t)}$ extends to a continuous map from $\bar{B}^4 \times [-\pi,\pi]$ into the space of surfaces (equipped with the flat topology).
\item For each point $a \in S^3$, there is a unique number $\tau(a)$ such that $\hat{\Sigma}_{(a,\tau(a))}$ is a totally geodesic two-sphere in $S^3$. For abbreviation, let $Q(a) \in \mathbb{RP}^3$ denote the unit normal vector to the surface $\hat{\Sigma}_{(a,\tau(a))}$.
\item If $\Sigma$ has genus at least $1$, then the map $Q: S^3 \to \mathbb{RP}^3$ has non-zero degree.
\end{itemize}
\end{proposition}

To see that the map $(a,t) \mapsto \hat{\Sigma}_{(a,t)}$ is continuous, one considers a sequence of points of the form $a_i = \Psi(x,r_i,s_i)$, where $x \in \Sigma$ and $\sqrt{r_i^2+s_i^2} \searrow \varepsilon$. Let $\tilde{a}_i = \Psi(x,\tilde{r}_i,\tilde{s}_i)$, where $\tilde{r}_i = (\frac{\sqrt{r_i^2+s_i^2}}{\varepsilon} - 1) \, r_i$ and $\tilde{s}_i = (\frac{\sqrt{r_i^2+s_i^2}}{\varepsilon} - 1) \, s_i$. As $i \to \infty$, the surfaces $\hat{\Sigma}_{(a_i,t)} = \Sigma_{(\tilde{a}_i,t)}$ converge to a geodesic sphere centered at the point $-\sin \theta \, x - \cos \theta \, \nu(x)$ of radius $\frac{\pi}{2} - \theta + t$, where $\theta \in [-\frac{\pi}{2},\frac{\pi}{2}]$ is defined by $\tan \theta = \lim_{i \to \infty} \frac{\tilde{s}_i}{\tilde{r}_i} = \lim_{i \to \infty} \frac{s_i}{r_i}$. In other words, $\sin \theta = \lim_{i \to \infty} \frac{s_i}{\sqrt{r_i^2+s_i^2}} = \lim_{i \to \infty} \frac{s_i}{\varepsilon}$. From this, the continuity property follows.

Finally, using the conformal invariance of the Willmore functional and a result of Ros \cite{Ros2} (cf. Proposition \ref{distance.surface} above), one obtains 
\[\sup_{(a,t) \in B^4 \times (-\pi,\pi)} \text{\rm area}(\Sigma_{(a,t)}) \leq \mathscr{W}(\Sigma),\] 
hence 
\begin{equation} 
\label{bound}
\sup_{(a,t) \in \bar{B}^4 \times [-\pi,\pi]} \text{\rm area}(\hat{\Sigma}_{(a,t)}) \leq \mathscr{W}(\Sigma). 
\end{equation}
We now state the main result in \cite{Marques-Neves}: 

\begin{theorem}[F.C.~Marques, A.~Neves \cite{Marques-Neves}]
If $\Sigma$ is an immersed minimal surface in $S^3$ of genus at least $1$, then $\text{\rm area}(\Sigma) \geq 2\pi^2$. Moreover, if $\Sigma$ is an arbitrary immersed surface in $S^3$ of genus at least $1$, then $\mathscr{W}(\Sigma) \geq 2\pi^2$.
\end{theorem}

The proof in \cite{Marques-Neves} is rather technical. In the following, we will merely sketch the main ideas. Suppose first that there exists an immersed minimal surface in $S^3$ which has genus at least $1$ and area less than $2\pi^2$. By Proposition \ref{lower.bound} above, any such surface must be embedded. Let $\mathscr{C}$ denote the set of all embedded minimal surfaces in $S^3$ which have genus at least $1$ and area less than $2\pi^2$. Clearly, $\mathscr{C} \neq \emptyset$. Moreover, it follows from results in \cite{Kuwert-Li-Schatzle} that the genus of a minimal surface in $\mathscr{C}$ is uniformly bounded from above (cf. Theorem \ref{high.genus} below). Using a theorem of Choi and Schoen \cite{Choi-Schoen}, we conclude that $\mathscr{C}$ is compact (see also Theorem \ref{choi.schoen.compactness} below). Consequently, there exists an embedded minimal surface $\Sigma \in \mathscr{C}$ which has smallest area among all surfaces in $\mathscr{C}$. 

Let $\mathscr{F}$ be the set of all continuous five-parameter families of surfaces $\{S_{(a,t)}: (a,t) \in \bar{B}^4 \times [-\pi,\pi]\}$ with the property that $S_{(a,t)} = \hat{\Sigma}_{(a,t)}$ for $(a,t) \in \partial(\bar{B}^4 \times [-\pi,\pi])$. Marques and Neves then define 
\begin{equation} 
\Lambda = \inf_{S \in \mathscr{F}} \sup_{(a,t) \in \bar{B}^4 \times [-\pi,\pi]} \text{\rm area}(S_{(a,t)}). 
\end{equation}
It is easy to see that $\Lambda \geq 4\pi$. Moreover, since $\Sigma$ is minimal, the inequality (\ref{bound}) gives 
\[\sup_{(a,t) \in \bar{B}^4 \times [-\pi,\pi]} \text{\rm area}(\hat{\Sigma}_{(a,t)}) \leq \text{\rm area}(\Sigma).\] 
Note that $\text{\rm area}(\Sigma) < 2\pi^2$, so $\Sigma$ cannot be congruent to the Clifford torus. Consequently, $\text{\rm ind}(\Sigma) \geq 6$ by Urbano's theorem. By perturbing the canonical family $\hat{\Sigma}_{(a,t)}$, one can construct a new five-parameter family of surfaces $S \in \mathscr{F}$ with the property that 
\[\sup_{(a,t) \in \bar{B}^4 \times [-\pi,\pi]} \text{\rm area}(S_{(a,t)}) < \text{\rm area}(\Sigma).\] 
Thus, 
\begin{equation} 
\label{ineq}
\Lambda < \text{\rm area}(\Sigma). 
\end{equation}
There are two cases now:

\textit{Case 1:} Suppose first that $\Lambda > 4\pi$. In this case, Marques and Neves show that there exists an embedded minimal surface $\tilde{\Sigma}$ with area $\text{\rm area}(\tilde{\Sigma}) = \Lambda$. In particular, $\text{\rm area}(\tilde{\Sigma}) = \Lambda > 4\pi$, so $\tilde{\Sigma}$ must have genus at least $1$. On the other hand, the inequality (\ref{ineq}) implies $\text{\rm area}(\tilde{\Sigma}) = \Lambda < \text{\rm area}(\Sigma)$. This contradicts the choice of $\Sigma$.

\textit{Case 2:} Suppose next that $\Lambda = 4\pi$. In this case, there exists a sequence $S^{(i)} \in \mathscr{F}$ such that 
\[\sup_{(a,t) \in \bar{B}^4 \times [-\pi,\pi]} \text{\rm area}(S_{(a,t)}^{(i)}) \leq 4\pi + \frac{1}{i}.\] 
For each $a \in B^4$ and each $i \in \mathbb{N}$, the surfaces $S_{(a,t)}^{(i)}$ form a sweepout of $S^3$. 

Let $V^{(i)}(a,t)$ denote the volume enclosed by the surface $S_{(a,t)}^{(i)}$, so that $V^{(i)}(a,-\pi) = 0$ and $V^{(i)}(a,\pi) = \text{\rm vol}(S^3)$. We may approximate the function $V^{(i)}(a,t)$ by a $C^1$-function $\tilde{V}^{(i)}(a,t)$ such that 
\[\sup_{(a,t) \in \bar{B}^4 \times [-\pi,\pi]} |V^{(i)}(a,t) - \tilde{V}^{(i)}(a,t)| \leq \frac{1}{i} \, \text{\rm vol}(S^3).\] 
By Sard's lemma we can assume that $\frac{1}{2} \, \text{\rm vol}(S^3)$ is a regular value of the function $(a,t) \mapsto \tilde{V}^{(i)}(a,t)$. 

Therefore, the set 
\[\Omega^{(i)} = \Big \{ (a,t) \in \bar{B}^4 \times [-\pi,\pi]: \tilde{V}^{(i)}(a,t) = \frac{1}{2} \, \text{\rm vol}(S^3) \Big \}\] 
is a smooth hypersurface in $\bar{B}^4 \times [-\pi,\pi]$. Moreover, one can arrange that the boundary $\partial \Omega^{(i)} \subset \partial B^4 \times (-\pi,\pi)$ is a graph over $\partial B^4$. In other words, for each $a \in B^4$ there is exactly one number $t \in (-\pi,\pi)$ such that $\tilde{V}^{(i)}(a,t) = \frac{1}{2} \, \text{\rm vol}(S^3)$.

Consider now a point $(a,t) \in \Omega^{(i)}$. Then the surface $S_{(a,t)}^{(i)}$ divides $S^3$ into two regions, each of which has volume at least $\big ( \frac{1}{2} - \frac{1}{i} \big ) \, \text{\rm vol}(S^3)$. On the other hand, we know that $\text{\rm area}(S_{(a,t)}^{(i)}) \leq 4\pi + \frac{1}{i}$. Hence, for each pair $(a,t) \in \Omega^{(i)}$, the surface $S_{(a,t)}^{(i)}$ is very close (in the sense of varifolds) to a totally geodesic sphere. Since the space of totally geodesic spheres is homeomorphic to $\mathbb{RP}^3$, one obtains a map $f^{(i)}: \Omega^{(i)} \to \mathbb{RP}^3$. Moreover, it turns out that $f^{(i)}(a,t) = Q(a)$ for each point $(a,t) \in (\partial B^4 \times [-\pi,\pi]) \cap \Omega_i$, where $Q$ is the map in Proposition \ref{canonical.family}. Hence, the map $Q: \partial B^4 \to \mathbb{RP}^3$ admits a continuous extension $f^{(i)}: \Omega^{(i)} \to \mathbb{RP}^3$. This contradicts the fact that $Q$ has non-zero degree. Therefore, any immersed minimal surface in $S^3$ of genus $1$ must have area at least $2\pi^2$. This proves the first statement.

The proof of the second statement involves similar ideas (see \cite{Marques-Neves} for details). Suppose that $\Sigma$ is an immersed surface of genus at least $1$ with Willmore energy less than $2\pi^2$. Proposition \ref{lower.bound} again implies that $\Sigma$ is embedded. Marques and Neves again define 
\[\Lambda = \inf_{S \in \mathscr{F}} \sup_{(a,t) \in \bar{B}^4 \times [-\pi,\pi]} \text{\rm area}(S_{(a,t)}).\] 
As above the inequality (\ref{bound}) gives 
\begin{equation} 
\label{ineq.2}
\Lambda \leq \sup_{(a,t) \in \bar{B}^4 \times [-\pi,\pi]} \text{\rm area}(\hat{\Sigma}_{(a,t)}) \leq \mathscr{W}(\Sigma) < 2\pi^2. 
\end{equation}
Moreover, the fact that $\Sigma$ has genus at least $1$ implies that $\Lambda > 4\pi$. Consequently, there exists an embedded minimal surface $\tilde{\Sigma}$ such that $\text{\rm area}(\tilde{\Sigma}) = \Lambda$. Since $\tilde{\rm area}(\tilde{\Sigma}) = \Lambda > 4\pi$, the surface $\tilde{\Sigma}$ must have genus at least $1$. On the other hand, it follows from (\ref{ineq.2}) that $\text{\rm area}(\tilde{\Sigma}) = \Lambda < 2\pi^2$. This contradicts the first statement. \\

We next mention a result concerning the Willmore energy of surfaces of high genus:

\begin{theorem}[E.~Kuwert, Y.~Li, R.~Sch\"atzle \cite{Kuwert-Li-Schatzle}]
\label{high.genus}
There exists a sequence of real numbers $\beta_g \in (4\pi,8\pi)$ such that $\lim_{g \to \infty} \beta_g = 8\pi$ and $\mathscr{W}(\Sigma) \geq \beta_g$ for every immersed surface $\Sigma$ in $S^3$ of genus $g$. In particular, $\text{\rm area}(\Sigma) \geq \beta_g$ for every immersed minimal surface $\Sigma$ in $S^3$ of genus $g$.
\end{theorem}

Finally, we note that Ilmanen and White \cite{Ilmanen-White} have recently obtained sharp estimates for the density of area-minimizing cones in Euclidean space. This result gives a lower bound for the area of certain minimal hypersurfaces in the unit sphere. 

\section{The first eigenvalue of the Laplacian on a minimal surface} 

In this final section, we describe an estimate for the first eigenvalue of the Laplace operator on a minimal surface. If $\Sigma$ is a minimal surface in $S^3$, then the restrictions of the coordinate functions in $\mathbb{R}^4$ satisfy 
\[\Delta_\Sigma x_i + 2 \, x_i = 0\] 
for $i \in \{1,2,3,4\}$. It was conjectured by Yau \cite{Yau} that the smallest positive eigenvalue of the operator $-\Delta_\Sigma$ is equal to $2$, provided that $\Sigma$ is embedded. While Yau's conjecture is an open problem, there are various partial results in this direction. In particular, the following result of Choi and Wang gives a lower bound for the first eigenvalue of the Laplacian on a minimal surface.

\begin{theorem}[H.I.~Choi and A.N.~Wang \cite{Choi-Wang}]
\label{eigenvalue.estimate}
Let $\Sigma$ be an embedded minimal surface in $S^3$, and let $\lambda$ be the smallest positive eigenvalue of the operator $-\Delta_\Sigma$. Then $\lambda > 1$. 
\end{theorem}

\textbf{Proof of Theorem \ref{eigenvalue.estimate}.} 
Suppose by contradiction that $\lambda \leq 1$. Let $\varphi: \Sigma \to \mathbb{R}$ be an eigenfunction, so that 
\[\Delta_\Sigma \varphi + \lambda \, \varphi = 0.\]
The surface $\Sigma$ divides $S^3$ into two regions, which we denote by $N$ and $\tilde{N}$. Let $\nu$ denote the outward-pointing unit normal vector field to $N$. Moreover, let $u: N \to \mathbb{R}$ and $\tilde{u}: \tilde{N} \to \mathbb{R}$ be harmonic functions satisfying $u|_\Sigma = \tilde{u}|_\Sigma = \varphi$. Using the Bochner formula, we obtain 
\[|D^2 u|^2 + 2 \, |\nabla u|^2 = \frac{1}{2} \, \Delta (|\nabla u|^2).\] 
We now integrating this identity over $N$ and apply the divergence theorem. This gives 
\begin{align*} 
&\int_N |D^2 u|^2 + 2 \int_N |\nabla u|^2 \\ 
&= \int_\Sigma \frac{1}{2} \, \langle \nabla(|\nabla u|^2),\nu \rangle \\ 
&= \int_\Sigma (D^2 u)(\nabla u,\nu) \\ 
&= \int_\Sigma \sum_{i=1}^2 (D^2 u)(e_i,\nu) \, \langle \nabla u,e_i \rangle + \int_\Sigma (D^2 u)(\nu,\nu) \, \langle \nabla u,\nu \rangle.  
\end{align*} 
Note that 
\[(D^2 u)(\nu,\nu) = -\sum_{i=1}^2 (D^2 u)(e_i,e_i) = -\Delta_\Sigma \varphi\] 
since $u$ is harmonic. We next define a function $\psi: \Sigma \to \mathbb{R}$ by $\psi = \langle \nabla u,\nu \rangle$. Then 
\[\langle \nabla^\Sigma \psi,e_i \rangle = (D^2 u)(e_i,\nu) + h(\nabla^\Sigma \varphi,e_i).\] 
Hence, we obtain 
\begin{align*} 
\int_N |D^2 u|^2 + 2 \int_N |\nabla u|^2 
&= \int_\Sigma \sum_{i=1}^2 (D^2 u)(e_i,\nu) \, \langle \nabla^\Sigma \varphi,e_i \rangle - \int_\Sigma \Delta_\Sigma \varphi \, \psi \\ 
&= \int_\Sigma \langle \nabla^\Sigma \varphi,\nabla^\Sigma \psi \rangle - \int_\Sigma h(\nabla^\Sigma \varphi,\nabla^\Sigma \varphi) - \int_\Sigma \Delta_\Sigma \varphi \, \psi \\ 
&= -\int_\Sigma h(\nabla^\Sigma \varphi,\nabla^\Sigma \varphi) - 2 \int_\Sigma \Delta_\Sigma \varphi \, \psi \\ 
&= -\int_\Sigma h(\nabla^\Sigma \varphi,\nabla^\Sigma \varphi) + 2\lambda \int_\Sigma \varphi \, \psi \\ 
&= -\int_\Sigma h(\nabla^\Sigma \varphi,\nabla^\Sigma \varphi) + 2\lambda \int_N |\nabla u|^2. 
\end{align*} 
Since $\lambda \leq 1$, we conclude that 
\[\int_N |D^2 u|^2 \leq -\int_\Sigma h(\nabla^\Sigma \varphi,\nabla^\Sigma \varphi).\] 
An analogous argument gives 
\[\int_{\tilde{N}} |D^2 \tilde{u}|^2 \leq \int_\Sigma h(\nabla^\Sigma \varphi,\nabla^\Sigma \varphi).\] 
(Note that the outward-pointing unit normal vector field to $\tilde{N}$ is given by $-\nu$, and the second fundamental form with respect to this choice of normal vector is $-h$.) Adding both identies gives 
\[\int_N |D^2 u|^2 + \int_{\tilde{N}} |D^2 \tilde{u}|^2 \leq 0.\] 
Therefore, $\nabla u$ is a parallel vector field on $N$. Substituting this back into the Bochner formula, we conclude that $\nabla u = 0$. Thus, $u$ is constant, and so is $\varphi$. This is a contradiction. \\

As a consequence of Theorem \ref{eigenvalue.estimate}, Choi and Schoen obtained a compactness theorem for embedded minimal surfaces in $S^3$.

\begin{theorem}[H.I.~Choi, R.~Schoen \cite{Choi-Schoen}]
\label{choi.schoen.compactness}
Given any integer $g \geq 1$, the space of all embedded minimal surfaces in $S^3$ of genus $g$ is compact.
\end{theorem}

Choe and Soret were recently able to verify Yau's conjecture for the Lawson surfaces and the Karcher-Pinkall-Sterling examples (cf. \cite{Choe}, \cite{Choe-Soret1}). The following result is a consequence of Courant's nodal theorem and plays an important role in the argument:

\begin{proposition}
\label{nodal}
Let $\Sigma$ be a closed surface equipped with a Riemannian metric. Let $\lambda$ be the smallest positive eigenvalue of the operator $-\Delta_\Sigma$, and let $\varphi$ be the associated eigenfunction. Moreover, let $\psi$ be another eigenfunction of the operator $-\Delta_\Sigma$ with eigenvalue $\mu > 0$. If $\{\psi=0\} \subset \{\varphi=0\}$, then $\lambda=\mu$.
\end{proposition}

\textbf{Proof of Proposition \ref{nodal}.} 
By assumption, we have $\{\varphi \neq 0\} \subset \{\psi \neq 0\}$. Hence, if we put $D_+ = \{\varphi > 0\}$, then we have 
\[D_+ = (D_+ \cap \{\psi > 0\}) \cup (D_+ \cap \{\psi < 0\}),\] 
Note that the sets $D_+ \cap \{\psi > 0\}$ and $D_+ \cap \{\psi < 0\}$ are disjoint open subsets of $\Sigma$. Moreover, the set $D_+$ is connected by Courant's nodal theorem (cf. \cite{Courant-Hilbert}, p.~452). Thus, we conclude that either $D_+ \cap \{\psi > 0\} = D_+$ or $D_+ \cap \{\psi < 0\} = D_+$. In other words, the restriction of $\psi$ to the set $D_+$ is either strictly positive or strictly negative. Similarly, we can show that the restriction of $\psi$ to the set $D_- = \{\varphi < 0\}$ is either strictly positive or strictly negative.

If $\psi|_{D_+}$ and $\psi|_{D_-}$ are of the same sign, then $\int_\Sigma \psi \neq 0$, which is impossible. Thus, $\psi|_{D_+}$ and $\psi|_{D_-}$ must have opposite signs. This implies that $\int_\Sigma \varphi \psi \neq 0$. Since $\varphi$ and $\psi$ are eigenfunctions of $-\Delta_\Sigma$ with eigenvalues $\lambda$ and $\mu$, we conclude that $\lambda=\mu$, as claimed. \\

\begin{corollary}[J.~Choe, M.~Soret \cite{Choe-Soret1}]
\label{invariance}
Let $\Sigma$ be an embedded minimal surface in $S^3$ which is symmetric under the reflection $\sigma(x) = x - 2 \, \langle a,x \rangle \, a$ for some unit vector $a \in \mathbb{R}^4$. Moreover, let $\lambda$ be the smallest positive eigenvalue of the operator $-\Delta_\Sigma$ and let $\varphi: \Sigma \to \mathbb{R}$ be an eigenfunction with eigenvalue $\lambda$. If $\lambda < 2$, then the eigenfunction $\varphi$ is invariant under the reflection $\sigma$.
\end{corollary}

\textbf{Proof of Corollary \ref{invariance}.} 
We argue by contradiction. Suppose that $\varphi \circ \sigma \neq \varphi$. Let $\tilde{\varphi} = \varphi \circ \sigma - \varphi$ and $\psi = \langle a,x \rangle$. Then $\tilde{\varphi}$ is an eigenfunction of $-\Delta_\Sigma$ with eigenvalue $\lambda$, and $\psi$ is an eigenfunction of $-\Delta_\Sigma$ with eigenvalue $2$. Moreover, we have $\{\psi=0\} \subset \{\tilde{\varphi} = 0\}$. Hence, Proposition \ref{nodal} implies that $\lambda=2$. This is a contradiction. \\

\begin{theorem}[J.~Choe, M.~Soret \cite{Choe-Soret1}]
\label{yau.conjecture.holds.for.known.surfaces}
Suppose that $\Sigma$ is one of the Lawson surfaces or one of the surfaces constructed by Karcher-Pinkall-Sterling. Then the smallest positive eigenvalue of $-\Delta_\Sigma$ is equal to $2$.
\end{theorem}

We will only give the proof of Theorem \ref{yau.conjecture.holds.for.known.surfaces} in the special case when $\Sigma$ is one of the Lawson surfaces. A key ingredient in the proof of Choe and Soret is the fact that the Lawson surfaces are invariant under reflection across certain geodesic two-spheres in $S^3$. To describe these symmetries, let us fix two positive integers $k$ and $m$. For each $i \in \mathbb{Z}_{2(k+1)}$, we consider the reflection $\sigma_i(x) = x - 2 \, \langle a_i,x \rangle \, a_i$, where 
\[a_i = \Big ( \sin \frac{\pi(2i+1)}{2(k+1)},-\cos \frac{\pi(2i+1)}{2(k+1)},0,0 \Big ).\] 
Similarly, for each $j \in \mathbb{Z}_{2(m+1)}$, we define $\tau_j(x) = x - 2 \, \langle b_j,x \rangle \, b_j$, where 
\[b_j = \Big ( 0,0,\sin \frac{\pi(2j+1)}{2(m+1)},-\cos \frac{\pi(2j+1)}{2(m+1)} \Big ).\] 
Let $\Gamma$ be the subgroup of $\text{\rm O}(4)$ generated by the reflections $\sigma_i$ and $\tau_j$. Note that the geodesic tetrahedron 
\begin{align*} 
T &= \Big \{ x \in S^3: -x_1 \, \sin \frac{\pi}{2(k+1)} < x_2 \, \cos \frac{\pi}{2(k+1)} \Big \} \\ 
&\cap \Big \{ x \in S^3: x_1 \, \sin \frac{\pi}{2(k+1)} > x_2 \, \cos \frac{\pi}{2(k+1)} \Big \} \\ 
&\cap \Big \{ x \in S^3: -x_3 \, \sin \frac{\pi}{2(m+1)} < x_4 \, \cos \frac{\pi }{2(m+1)} \Big \} \\ 
&\cap \Big \{ x \in S^3: x_3 \, \sin \frac{\pi}{2(m+1)} > x_4 \, \cos \frac{\pi}{2(m+1)} \Big \} 
\end{align*} 
is a fundamental domain for $\Gamma$.

As in the proof of Theorem \ref{lawson.surfaces}, let $\Sigma_{0,0}$ be an embedded least area disk whose boundary is the geodesic quadrilateral with vertices $P_0$, $Q_0$, $P_1$, and $Q_1$. Note that the reflections $\sigma_0$ and $\tau_0$ map the boundary $\partial \Sigma_{0,0}$ to itself. 

\begin{lemma}
\label{symmetry.of.fundamental.piece}
The surface $\Sigma_{0,0}$ is invariant under the reflections $\sigma_0$ and $\tau_0$.
\end{lemma} 

\textbf{Proof of Lemma \ref{symmetry.of.fundamental.piece}.} 
It suffices to show that $\Sigma_{0,0}$ is invariant under the reflection $\sigma_0$. Note that the set $\{x \in \Sigma_{0,0}: \langle a_0,x \rangle = 0\}$ is a union of finitely many smooth arcs. In the first step, we show that the set $\Sigma_+ = \{x \in \Sigma_{0,0}: \langle a_0,x \rangle > 0\}$ is connected. Indeed, if $\Sigma_+$ is disconnected, then there exists a connected component of $\Sigma_+$ which is disjoint from the boundary $\partial \Sigma_{0,0}$. Let us denote this connected component by $D$. Clearly, $D$ is a stable minimal surface whose boundary is contained in the totally geodesic two-sphere $\{x \in S^3: \langle a_0,x \rangle = 0\}$. Using the function $\langle a_0,x \rangle$ as a test function in the stability inequality, we conclude that the second fundamental form vanishes at each point on $D$ (see \cite{Ros1}, Lemma 1). Therefore, $\Sigma_{0,0}$ is totally geodesic, which is impossible. Thus, $\Sigma_+$ is connected. An analogous argument shows that the set $\Sigma_- = \{x \in \Sigma_{0,0}: \langle a_0,x \rangle < 0\}$ is connected as well. Since $\Sigma_{0,0}$ is homeomorphic to a disk, we conclude that $\Sigma_+$ and $\Sigma_-$ are simply connected.

After replacing $\Sigma_{0,0}$ by $\sigma_0(\Sigma_{0,0})$ if necessary, we can arrange that $\text{\rm area}(\Sigma_+) \leq \frac{1}{2} \, \text{\rm area}(\Sigma_{0,0})$. The surface $\bar{\Sigma}_+ \cup \sigma_0(\bar{\Sigma}_+)$ is homemorphic to a disk, and its area is bounded from above by the area of $\Sigma_{0,0}$. Consequently, the surface $\bar{\Sigma}_+ \cup \sigma_0(\bar{\Sigma}_+)$ is a least area disk. In particular, the surface $\bar{\Sigma}_+ \cup \sigma_0(\bar{\Sigma}_+)$ is smooth and has zero mean curvature. Hence, the unique continuation theorem implies that $\bar{\Sigma}_+ \cup \sigma_0(\bar{\Sigma}_+) = \Sigma_{0,0}$. This shows that $\Sigma_{0,0}$ is invariant under the reflection $\sigma_0$, thus completing the proof of Lemma \ref{symmetry.of.fundamental.piece}. \\

After these preparations, we now describe the proof of Theorem \ref{yau.conjecture.holds.for.known.surfaces}. Let $\Sigma = \bigcup_{(i,j) \in A_{\text{\rm even}}} \Sigma_{i,j}$ be the Lawson surface constructed in Theorem \ref{lawson.surfaces}. It follows from Lemma \ref{symmetry.of.fundamental.piece} that the surface $\Sigma$ is invariant under $\Gamma$. Let $F: B^2 \to \Sigma_{0,0}$ be a conformal parametrization of $\Sigma_{0,0}$. After composing $F$ with a M\"obius transformation on $B^2$, we can arrange that $F^{-1} \circ \sigma_0 \circ F(s,t) = (-s,t)$ and $F^{-1} \circ \tau_0 \circ F(s,t) = (s,-t)$. Thus, the pre-image of the surface $\Sigma_{0,0} \cap \overline{T}$ under the map $F$ is a quadrant in $B^2$. From this, we deduce that the fundamental patch 
\[S = \Sigma \cap \overline{T} = (\Sigma_{0,0} \cap \overline{T}) \cup (\Sigma_{-1,-1} \cap \overline{T})\] 
is simply connected, and the intersection of $S$ with each face of $T$ is a connected curve.

Let $\lambda$ be the smallest positive eigenvalue of the operator $-\Delta_\Sigma$, and let $\varphi$ be an associated eigenfunction. If $\lambda < 2$, then $\varphi$ is invariant under $\Gamma$ by Corollary \ref{invariance}. The nodal set $\{\varphi=0\}$ is a union of finitely many smooth arcs. Let us choose a piecewise smooth curve $C \subset S \cap \{\varphi=0\}$ which starts at a point on the boundary $\partial S$ and ends at another point on the boundary $\partial S$. There exists a connected component of $S \setminus C$ which is disjoint from one of the faces of $T$. Let us denote this connected component by $D$, and let $D'$ be another connected component of $S \setminus D$ which is disjoint from $D$. By assumption, we have $D \cap F = \emptyset$, where $F$ is one of the faces of the geodesic tetrahedron $T$. 

Let us pick two points $x$ and $y$ in the interior of $S$ such that $x \in D \cap \{\varphi \neq 0\}$ and $y \in D' \cap \{\varphi \neq 0\}$. Finally, let $z \in \{\varphi \neq 0\}$ denote the reflection of $x$ across $F$. By Courant's nodal theorem, two of the points $x,y,z$ lie in the same connected component of $\{\varphi \neq 0\}$. There are three cases now:

\textit{Case 1:} Suppose that $x$ and $y$ lie in the same connected component of $\{\varphi \neq 0\}$. Let $\alpha: [0,1] \to \{\varphi \neq 0\}$ be a continuous path such that $\alpha(0) = x$ and $\alpha(1) = y$. We can find a continuous path $\tilde{\alpha}: [0,1] \to S \cap \{\varphi \neq 0\}$ with the property that $\tilde{\alpha}(0) = x$ and $\tilde{\alpha}(t) = \rho(t) \, \alpha(t)$ for some element $\rho(t) \in \Gamma$. Clearly, $\tilde{\alpha}(1) = y$. Since the path $\tilde{\alpha}(t)$ cannot intersect $C$, it follows that $x$ and $y$ belong to the same connected component of $S \setminus C$. This contradicts our choice of $x$ and $y$.

\textit{Case 2:} Suppose that $y$ and $z$ lie in the same connected component of $\{\varphi \neq 0\}$. Let $\alpha: [0,1] \to \{\varphi \neq 0\}$ be a continuous path such that $\alpha(0) = z$ and $\alpha(1) = y$. In this case, there exists a continuous path $\tilde{\alpha}: [0,1] \to S \cap \{\varphi \neq 0\}$ such that $\tilde{\alpha}(0) = x$ and $\tilde{\alpha}(t) = \rho(t) \, \alpha(t)$ for some element $\rho(t) \in \Gamma$. Clearly, $\tilde{\alpha}(1) = y$. Since the path $\tilde{\alpha}(t)$ cannot intersect $C$, it follows that $x$ and $y$ belong to the same connected component of $S \setminus C$. This contradicts our choice of $x$ and $y$.

\textit{Case 3:} Suppose that $x$ and $z$ lie in the same connected component of $\{\varphi \neq 0\}$. Let $\alpha: [0,1] \to \{\varphi \neq 0\}$ be a continuous path such that $\alpha(0) = x$ and $\alpha(1) = z$. We can find a continuous path $\tilde{\alpha}: [0,1] \to S \cap \{\varphi \neq 0\}$ with the property that $\tilde{\alpha}(0) = x$ and $\tilde{\alpha}(t) = \rho(t) \, \alpha(t)$ for some element $\rho(t) \in \Gamma$. Clearly, $\tilde{\alpha}(1) = x$. Moreover, since the path $\tilde{\alpha}(t)$ cannot intersect $C$, we conclude that the path $\tilde{\alpha}(t)$ is disjoint from $F$. From this, we deduce that $\rho(t) \in \Gamma_0$, where $\Gamma_0$ denotes the subgroup of $\Gamma$ which is generated by the reflections across the faces of $T$ different from $F$. On the other hand, the identity $x = \tilde{\alpha}(1) = \rho(1) \, \alpha(1) = \rho(1) \, z$ implies $\rho(1) \notin \Gamma_0$. Again, this is a contradiction. This completes the proof of Theorem \ref{yau.conjecture.holds.for.known.surfaces}. \\

Finally, let us mention the following theorem due to Ros \cite{Ros1}: 

\begin{theorem}[A.~Ros \cite{Ros1}]
\label{two.piece}
Let $\Sigma$ be an embedded minimal surface in $S^3$, and let $a$ be a unit vector in $\mathbb{R}^4$. Then the set $\{x \in \Sigma: \langle a,x \rangle > 0\}$ is connected.
\end{theorem}

\textbf{Sketch of the proof of Theorem \ref{two.piece}.} 
If $\Sigma$ is a totally geodesic two-sphere, the assertion is trivial. We will, therefore, assume that $\Sigma$ is not totally geodesic. Let $D$ be a connected component of the set $\{x \in \Sigma: \langle a,x \rangle > 0\}$, and let $\Gamma$ denote the boundary of $D$, so that $\Gamma \subset \{x \in S^3: \langle a,x \rangle = 0\}$. Since $D$ is a nodal domain of an eigenfunction of the Laplace operator, the boundary $\Gamma$ is a union of finitely many smooth arcs. The surface $\Sigma$ divides $S^3$ into two regions, which we denote by $N$ and $\tilde{N}$. Note that the regions $N$ and $\tilde{N}$ are mean convex. Since the curve $\Gamma$ is null-homologous in $N$, we can find an area-minimizing surface $S \subset N$ such that $\partial S = \Gamma$. Similarly, there exists an area-minimizing surface $\tilde{S} \subset \tilde{N}$ satisfying $\partial \tilde{S} = \Gamma$. Note that $\Gamma$ may not be connected, and $S$ and $\tilde{S}$ might be disconnected as well. Using the function $\langle a,x \rangle$ as a test function in the stability inequality, we conclude that $S$ is totally geodesic (see \cite{Ros1}, Lemma 1). An analogous argument shows that $\tilde{S}$ is totally geodesic. 

We now distinguish two cases: 

\textit{Case 1:} Suppose first that 
\[S \cup \tilde{S} \subset \{x \in S^3: \langle a,x \rangle = 0\}.\] 
Since $\partial S = \partial \tilde{S} = \Gamma$, we have 
\[S \cup \tilde{S} = \{x \in S^3: \langle a,x \rangle = 0\}.\] 
Since the surfaces $S$ and $\tilde{S}$ cannot touch $\Sigma$ away from $\Gamma$, we conclude that 
\[\{x \in \Sigma: \langle a,x \rangle = 0\} = (S \cup \tilde{S}) \cap \Sigma = \Gamma.\] 
This shows that $\{x \in \Sigma: \langle a,x \rangle > 0\} = D$, as claimed.

\textit{Case 2:} Suppose finally that 
\[S \cup \tilde{S} \not\subset \{x \in S^3: \langle a,x \rangle = 0\}.\] 
Without loss of generality, we may assume that 
\[S \not\subset \{x \in S^3: \langle a,x \rangle = 0\}.\] 
Let $S_0$ be connected component of $S$ such that 
\[S_0 \not\subset \{x \in S^3: \langle a,x \rangle = 0\}.\] 
Since $S_0$ is totally geodesic, we have 
\[S_0 \subset \{x \in S^3: \langle b,x \rangle = 0\}\] 
for some unit vector $b \neq a$. This implies 
\[\partial S_0 \subset \partial S \cap \{x \in S^3: \langle b,x \rangle = 0\} \subset \{x \in S^3: \langle a,x \rangle = \langle b,x \rangle = 0\}.\] 
Thus, $S_0$ is a totally geodesic hemisphere. Moreover, $S_0$ does not touch $\Sigma$ except along the boundary. We now rotate the surface $S_0$ until it touches $\Sigma$. When that happens, the two surfaces conincide by the strict maximum principle. In particular, it follows that $\Sigma$ is totally geodesic, contrary to our assumption. This completes the proof of Theorem \ref{two.piece}. \\

Note that, if Yau's conjecture is true, then the function $\langle a,x \rangle$ is a first eigenfunction of the operator $-\Delta_\Sigma$, and Theorem \ref{two.piece} is a consequence of Courant's nodal theorem. 

We remark that many results in this section can be extended to higher dimensions. For example, the eigenvalue estimate of Choi and Wang works in all dimensions. Moreover, the two-piece property was generalized to higher dimensions in \cite{Choe-Soret1}. Finally, Tang and Yan \cite{Tang-Yan} recently obtained a sharp eigenvalue estimate for isoparametric minimal surfaces in $S^n$.

\end{document}